\documentclass[10pt]{extarticle}

\usepackage[margin=1.5in]{geometry}
\usepackage{amsthm,amsmath,amssymb,amsfonts,latexsym,enumitem,hyperref,stackengine,tikz-cd,titlesec}% more math symbols not covered by amssymb
\usepackage[autostyle]{csquotes}% for quoting things aesthetically
\usepackage[style=alphabetic, sorting=ynt, maxbibnames=99, backend=biber]{biblatex}

\titleformat{\section}[hang]{\normalfont\large\bfseries\centering}{§\,\thesection}{1em}{}
\titleformat{\subsection}[hang]{\normalfont\bfseries}{§\,\thesubsection}{1em}{}
\hypersetup{colorlinks,citecolor=black,filecolor=black,linkcolor=black,urlcolor=black}% for styling links
\newcommand\blfootnote[1]{%
	\begingroup
	\renewcommand\thefootnote{}\footnote{#1}%
	\addtocounter{footnote}{-1}%
	\endgroup
}

\newcommand{\ptsabove}{1em}
\newcommand{\ptsbelow}{1em}
\newcommand{\ptsindentamount}{0em}
\newcommand{\theoremmarin}{0em}
\newcommand{\emsbelowtheoremhead}{0.5em}

\newtheoremstyle{default theorem style}{\ptsabove}{\ptsbelow}{\addtolength{\leftskip}{\theoremmarin}\itshape}{\ptsindentamount}{\bf}{.}{\emsbelowtheoremhead}{(\thmnumber{#2}) \thmname{#1}}
\newtheoremstyle{theorem style for definitions}{\ptsabove}{\ptsbelow}{}{\ptsindentamount}{\bf}{.}{\emsbelowtheoremhead}{(\thmnumber{#2}) \thmname{#1}}

\theoremstyle{default theorem style}
\newtheorem{theorem}{Theorem}[section]
\newtheorem{lemma}[theorem]{Lemma}
\newtheorem{claim}[theorem]{Claim}
\newtheorem{prop}[theorem]{Proposition}
\newtheorem{cor}{Corollary}[theorem]

\theoremstyle{theorem style for definitions}

\newtheorem{example}[theorem]{Example}
\newtheorem{examples}[theorem]{Examples}

\newtheorem{conjecture}[theorem]{Conjecture}

\setlist[enumerate,1]{label={(\roman*)}}% for formating enumerated lists.

\newcommand{\Z}{{\mathbb Z}}
\newcommand{\ass}{\mathrm{ass}}
\newcommand{\m}{\mathfrak{m}}% for maximal ideals
\newcommand{\p}{\mathfrak{p}}% for prime ideals
\newcommand{\q}{\mathfrak{q}}% for prime ideals
\newcommand{\seq}{\subseteq}
\renewcommand{\>}{\rangle}
\newcommand{\<}{\langle}
\renewcommand{\tilde}{\widetilde}% so tilde is more pronounced

\title{\textbf{\large{Collections of an Ideal:\\ Any $n$-Absorbing Ideal is Strongly $n$-Absorbing}}}
\date{\vspace{-4ex}}
\author{Spencer Secord}

\begin{document}
	\maketitle
	\begin{abstract}
		We show that any $n$-absorbing ideal must be strongly $n$-absorbing, which is the first of Anderson and Badawi's three interconnected conjectures on absorbing ideals \cite{n-abs}. 
		We prove this by introducing and studying objects called maximal and semimaximal collections, which are tools for analysing the multiplicative ideal structure of commutative rings. 
		We also make heavy use of previous work on absorbing ideals done by Anderson--Badawi \cite{n-abs}, Choi--Walker \cite{Choi_Walker}, Laradji \cite{Laradji}, and Donadze \cite{Donadze_radical_formula}, among others.
		At the end of this paper we pose three conjectures, each of which imply the last unsolved conjecture made by Anderson and Badawi.
	\end{abstract}
	
	
	\blfootnote{\textit{Email:} spencer.e.secord+math@gmail.com}
	\vspace{-3em}
	
	\section{Introduction}\label{section: intro}
	For an ideal $\alpha$ of a commutative ring $A$, we define a \textit{collection} $C$ of $\alpha$ to be a list $C=[C_1,...,C_n]$ containing ideals $C_i$ of $A$ such that $\prod_i C_i\seq \alpha$ but $\prod_{i\ne j} C_i\not\seq \alpha$ for all $j$. We say $C$ is an \textit{$n$-collection} whenever $C=[C_1,...,C_n]$ is a collection whose size $|C|$ is equal to $n$.
	
	We will primarily focus our attention on two specific types of collections: maximal and semimaximal collections.
	Given collections $C$ and $C'$ of $\alpha$ with $|C|\leq |C'|$, we will say $C$ \textit{lies below} $C'$, or $C\leq C'$, whenever $C'$ can be reordered such that $C_i\seq C'_i$ for all $C_i$. 
	A \textit{maximal} collection $M$ of $\alpha$ is a collection which is maximal amongst all collections of $\alpha$ (up to reordering). Equivalently, a maximal collections $M$ is any collection which only lies below collections which are equivalent to $M$ after reordering. 
	Similarly, a \textit{semimaximal} collection $N$ is a collection which is maximal amongst collections of the same size.
	In other words, the only collections of the same size lying above it are reorderings of itself.
	
	For example, consider the ideal $\alpha=\<x,y,z\>^2\cap\<x^2,y\>$ of the ring $\mathbb{Q}[x,y,z]$, which is discussed in (\ref{example: classifying maximal collections}). Up to reordering, $\alpha$ contains two maximal collections: $[\<x^2,y\>,\<x,y,z\>]$ and $[\<x,y\>,\<x,y\>]$. The only semimaximal collection of $\alpha$ which is not maximal is the trivial semimaximal collection $[\alpha]$.
	
	Given any collection $C$ of $\alpha$, the list $C_{i\ne j}\cup[(\alpha:\prod_{i\ne j}C_i)]$ will always be a collection of the same size which lies above $C$, where \enquote{$\cup$} represents concatenation and $C_{i\ne j}$ is the sublist of $C$ which misses the $j^{th}$ element $C_j$. 
	Therefore, a collection $N$ is semimaximal if and only if $N_j = (\alpha:\prod_{i\ne j}N_i)$ for every element $N_j$ of $N$.
	
	In this paper, we will be using results on collections to study $n$-absorbing and strongly $n$-absorbing ideals, which will be defined shortly. 
	In \cite{n-abs}, Anderson and Badawi showed that every ideal of a Noetherian ring will be (strongly) $n$-absorbing for some positive integer $n$, 
	and thus the theory of absorbing ideals is applicable to any ideal of a Noetherian ring.
	
	Absorbing ideals were first introduced by Badawi \cite{2-abs} as a generalization of prime ideals: an \textit{$n$-absorbing ideal} is any ideal whose principal collections $[\<f_i\>]_i$ are of size at most $n$.
	Equivalently, $\alpha$ is \textit{n-absorbing} if for every $f_1,...,f_{n+1}\in A$ with $\prod_i f_i \in \alpha$, there exists some $j$ such that $\prod_{i\ne j}f_i\in \alpha$. 
	Hence $\alpha$ is $1$-absorbing if and only if $\alpha$ is a prime ideal.
	We use $\omega_A(\alpha)$ to represent the smallest integer $n$ such that $\alpha$ is an $n$-absorbing ideal of $A$, and we say $\omega_A(\alpha)=\infty$ if no such $n$ exists. We will write $\omega(\alpha)$ when no confusion about $A$ can be made.
	
	Later, Anderson and Badawi \cite{n-abs} studied absorbing ideals in more depth, and introduced the notion of a strongly $n$-absorbing ideal: an ideal $\alpha$ is \textit{strongly $n$-absorbing} if every collection of $\alpha$ has size at most $n$.
	Equivalently, $\alpha$ is strongly n-absorbing if for every list of $n+1$ ideals $[\gamma_i]_{i=1}^{n+1}$ satisfying $\prod_i \gamma_i \seq \alpha$, there exists some $j$ such that $\prod_{i\ne j}\gamma_i\seq \alpha$. 
	We write $\omega_A^*(\alpha)$ to represent the smallest integer $n$ such that $\alpha$ is strongly $n$-absorbing, and we set $\omega_A^*(\alpha)=\infty$ if no such $n$ exists.
	We call $\alpha$ an \textit{absorbing ideal} whenever $\omega^*(\alpha)<\infty$.
	
	The main result of this paper is a proof that these two concepts are equivalent; that every $n$-absorbing ideal is strongly $n$-absorbing, and thus $\omega(\alpha)=\omega^*(\alpha)$ for any ideal $\alpha$.
	Anderson and Badawi hypothesised this in \cite{n-abs}, their first of three interconnected conjectures which gained further attention after the conjectures appeared in \cite[Problem 30]{open}. 
	Only a solution to the third conjecture has been published so far, first \cite{Choi_Walker_arxiv} by Choi and Walker \cite{Choi_Walker}, but also independently by Donadze \cite{Donadze_radical_formula}. 
	We state each conjecture below, numbered as in \cite{Badawi-overview-of-conjectures}.
	\begin{enumerate}[label={}, leftmargin = 10pt]
		\item \hypertarget{(C1)}{\textbf{(C1).}} Any n-absorbing ideal is strongly n-absorbing.
		\item \hypertarget{(C2)}{\textbf{(C2).}} If $\alpha$ is an n-absorbing ideal of $A$, then $\omega_{A[x]}(\alpha[x])=\omega(\alpha)$.
		\item \hypertarget{(C3)}{\textbf{(C3).}} If $\alpha$ is n-absorbing, then $\left(\sqrt{\alpha}\right)^n \seq\alpha$.
	\end{enumerate}
	A grossly oversimplified summary of the current progress on the three conjectures can be visualized with the following diagram:
	
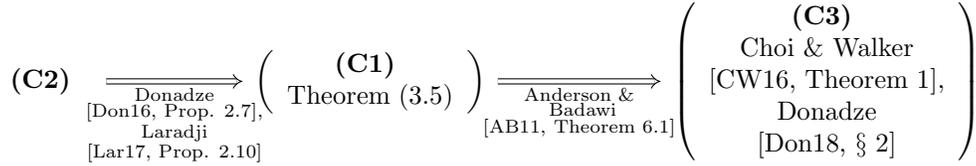
\vspace{1.75em}
	\begin{center}
		\begin{tikzcd}[remember picture, overlay]
			\begin{array}{c}
				\textup{\textbf{\hyperlink{(C2)}{(C2)} }}
			\end{array}
			\arrow[rrr,Rightarrow, "\substack{\substack{\textup{Donadze}\\\substack{\textup{\cite[Prop. 2.7]{Donadze_infinite_fields},}\\\substack{\textup{Laradji}\\\textup{\cite[Prop. 2.10]{Laradji}}}}} }
			"']
			&[0em] &[-2em] &[0em]
			\left(
			\begin{array}{c}
				\textup{\textbf{\hyperlink{(C1)}{(C1)} }}\\
				\textup{Theorem (\ref{Theorem: (C1) - the First Anderson-Badawi Conjecture})}
			\end{array}
			\right)
			\arrow[rrr,Rightarrow,"\substack{\textup{Anderson $\&$}\\\substack{\textup{ Badawi}\\\textup{\cite[Theorem 6.1]{n-abs}}}}"']  &[0em] &[-1em] &[0em]
			\left(
			\begin{array}{c}
				\textup{\textbf{\hyperlink{(C3)}{(C3)} }}\\
				\textup{Choi $\&$ Walker}\\
				\textup{\cite[Theorem 1]{Choi_Walker_arxiv},}\\
				\textup{Donadze}\\
				\textup{\cite[$\S$ 2]{Donadze_radical_formula}}
			\end{array}\right)
		\end{tikzcd}
	\end{center}
	\vspace{3em}
	For an in-depth review of the previous methods and progress on the three conjectures, including references to proofs of (C1) in several specific cases, we refer the reader to \cite{Badawi-overview-of-conjectures}. 
	Before continuing, we include the table of contents for the reader's convenience.
	\tableofcontents

	\subsection{Summary}\label{subsection: summary}
	Section \ref{section: intro} aims to provide the reader a brief introduction to collections, absorbing ideals, and the three conjectures given by Anderson and Badawi. 
	In particular, Subsection \ref{section: preliminaries on absorbing ideals} gives a self-contained treatment of some of the strongest and most fundamental results on $n$-absorbing ideals appearing in the literature. 
	
	Section \ref{section: collections of an ideal} introduces the reader to some basic techniques relating to collections, with particular attention given to maximal and semimaximal collections. 
	In Subsection \ref{subsection: ontro to maximal collections} we discuss general results on maximal and semimaximal collections, while giving examples and motivation.
	Subsection \ref{subsection: tilde ideals} introduces the tilde ideals $\tilde{C_I}$, which are a useful tool to generalize results on maximal and semimaximal collections.
	The new techniques made available by the $\tilde{C_I}$ constructions will prove useful in the last two sections of this paper.
	In Subsection \ref{subsection: artin rings} we restrict to the case of Artin rings, and prove that there exists a unique (up to reordering) maximal collection of the zero ideal. 
	This result gives us a kind of \enquote{base case} for our proof of the main result of this paper: Theorem (\ref{Theorem: (C1) - the First Anderson-Badawi Conjecture}).
	
	The purpose of Section \ref{section: the first conjecture} is to prove our main result, Theorem (\ref{Theorem: (C1) - the First Anderson-Badawi Conjecture}), which states that any $n$-absorbing ideal is strongly $n$-absorbing. 
	This is the statement of Anderson and Badawi's first conjecture \hyperlink{(C1)}{(C1)}.
	Lemma (\ref{Lemma: reduction to quasilocal noetherian rings}) is a key step in our approach, which allows us to restrict our focus to strongly $n$-absorbing ideals $\alpha$ of a quasilocal Noetherian ring whose maximal $n$-collections must contain a copy of each maximal ideal. 
	After restricting to this context, we prove Theorem (\ref{Theorem: (C1) - the First Anderson-Badawi Conjecture}) by constructing an $n$-collection of principal ideals for $\alpha$.
	
	Section \ref{section: some new conjectures} discusses three new conjectures which are stronger or equivalent to the second \hyperlink{(C2)}{(C2)} and only remaining conjecture by Anderson and Badawi.
	Subsection \ref{subsection: collections of alpha[x]} starts by exploring the strictness of the inclusion $\tilde C_j\seq c(\tilde C_j)[X]$, culminating in the introduction of Conjecture (\ref{conjecture: maximal collections of alpha[x] are M[x]}). 
	This conjecture hypothesizes that for any ideal $\alpha$ of a Noetherian ring, the maximal collections of $\alpha[x]$ are of the form $[M_i[x]]_i$ for some maximal collection $[M_i]_i$ of $\alpha$. 
	In the following discussion, we show that Conjecture (\ref{conjecture: maximal collections of alpha[x] are M[x]}) is stronger than \hyperlink{(C2)}{(C2)}.
	Subsection \ref{subsection: finitely many maximal collections} focuses on Conjecture (\ref{conjecture: finitely many maximal collections}), which states that any ideal of a Noetherian ring will have finitely many maximal collections.
	Later in this subsection, we show that this conjecture is stronger than (\ref{conjecture: maximal collections of alpha[x] are M[x]}), and therefore must imply \hyperlink{(C2)}{(C2)}.
	This is surprising since Conjecture (\ref{conjecture: finitely many maximal collections}) was formulated with the motivation of exploring the feasibility of a computational study of maximal collections and absorbing ideals. 
	Finally, Subsection \ref{subsection: absorbing ideals in graded rings} focuses on our last Conjecture (\ref{conjecture: graded ring formulation of Anderson-Badawi}), which postulates that if $\alpha$ is a $n$-absorbing ideal of a graded ring $R$, then the largest homogeneous ideal $\alpha^*$ contained in $\alpha$ must also be an $n$-absorbing ideal. This is a natural generalization of the well known result: if $\p$ is a prime ideal in a graded ring $R$, then $\p^*$ is also a prime ideal of $R$.
	The rest of this subsection is dedicated to proving that Conjecture (\ref{conjecture: graded ring formulation of Anderson-Badawi}) is equivalent to \hyperlink{(C2)}{(C2)}.
	The main results of Section \ref{section: some new conjectures} are summarized in the following diagram:
	\begin{center}
		\begin{tikzcd}[remember picture, overlay]
			\begin{array}{c}
				\textup{\text{(\ref{conjecture: finitely many maximal collections})}}
			\end{array}
			\arrow[rrr,Rightarrow, "\text{Prop. } (\ref{prop: conjecture finitely many maximal collections implies conjecture maximal collections of alpha[x] are M[x]})
			"']
			&[0em] &[-2em] &[0em]
			\begin{array}{c}
				\textup{\text{(\ref{conjecture: maximal collections of alpha[x] are M[x]})}}
			\end{array}
			\arrow[r,Rightarrow, ""']&
			\begin{array}{c}
				\textup{\textbf{\hyperlink{(C2)}{(C2)}}}
			\end{array}
			\arrow[rrr,Leftrightarrow,"\text{Cor. }(\ref{question two is equivalent to the Anderson-Badawi conjectures})"']  &[0em] &[-1.5em] &[0em]%
			\begin{array}{c}
				\textup{\text{(\ref{conjecture: graded ring formulation of Anderson-Badawi})}}
			\end{array}
		\end{tikzcd}
	\end{center}

	\subsection{Preliminaries on Absorbing Ideals}\label{section: preliminaries on absorbing ideals}
	Let $\alpha$ be an ideal of a ring $A$. 
	Clearly we will always have $\omega(\alpha)\leq\omega^*(\alpha)$.
	Given a subideal $\beta\seq \alpha$, we have $\omega(\alpha)=\omega(\alpha/\beta)$ and $\omega^*(\alpha)=\omega^*(\alpha/\beta)$. For any multiplicatively closed subset $S$ of $A$ we have $\omega(S^{-1}\alpha)\leq \omega(\alpha)$ and $\omega^*(S^{-1}\alpha)\leq \omega^*(\alpha)$. 
	If $\gamma_1$ and $\gamma_2$ are ideals of $A$, then $\omega(\gamma_1\cap \gamma_2)\leq \omega(\gamma_1)+\omega(\gamma_2)$ and $\omega^*(\gamma_1\cap \gamma_2)\leq \omega^*(\gamma_1)+\omega^*(\gamma_2)$. 
	Therefore $\omega(\bigcap_i \gamma_i)\leq \sum_i \omega(\gamma_i)$ and $\omega^*(\bigcap_i \gamma_i)\leq \sum_i \omega^*(\gamma_i)$ for any family of ideals $\{\gamma_i\}_i$ of $A$.\par
	An ideal $\alpha$ is strongly $n$-absorbing if and only if $\alpha$ is $n$-absorbing for finitely generated ideals; 
	i.e. $\omega^*(\alpha)\leq n$ if and only if $\alpha$ has no collections of size greater than $n$ whose elements are finitely generated ideals. 
	Indeed, if $C$ is an $n$-collection of $\alpha$, then for all $j$ we have $\prod_{i\ne j} C_i\not\seq \alpha$, and thus for every $i\ne j$ there must exist $f_{ij}\in C_i$ such that $\prod_{i\ne j}f_{ij}\not\in \alpha$. 
	Hence $\prod_i \<\{f_{ij}\}_{j=1,j\ne i}^n\> \seq \alpha$ but $\prod_{i\ne j} \<\{f_{ij}\}_{j=1,j\ne i}^n\> \not\seq \alpha$ for all $j$. 
	It follows that $[\<\{f_{ij}\}_{j\ne i}\>]_i$ is an $n$-collection of $\alpha$ containing finitely generated ideals.
	
	If $B$ is a subring of $A$, then $\omega_B(\alpha\cap B)\leq\omega(\alpha)$ and $\omega_B^*(\alpha\cap B)\leq\omega^*(\alpha)$. 
	Let $B'$ be the subring of $A$ generated by the finite set $\{f_{ij}\}_{ij}=\bigcup_j \{f_{ij}\}_{i\ne j}$ defined in the previous paragraph. 
	Clearly $[\<\{f_{ij}\}_{j\ne i}\>_{B'}]_i$ is an $n$-collection of $\alpha\cap B'$, and so $n\leq\omega_{B'}^*(\alpha\cap B')$. 
	It follows logically that $\alpha$ is strongly $n$-absorbing if and only if $\alpha\cap B$ is a strongly $n$-absorbing ideal of $B$ for every finitely generated subring $B$ of $A$. 
	Similarly, $\alpha$ is $n$-absorbing if and only if $\alpha\cap B$ is an $n$-absorbing ideal of $B$ for every finitely generated subring $B$ of $A$. 
	
	In the later parts of this paper, we will often study an ideal $\alpha$ of some ring $A$ by looking at ideals $\alpha\cap B$ for appropriately large finitely generated subrings $B$ of $A$. 
	The main advantage of working with finitely generated subrings is that they are isomorphic to a quotient of $\Z[x_1,...,x_n]$, and will therefore be Noetherian. 
	The Noetherian condition is particularly useful when studying maximal collections, in no small part due to the connection between maximal collections and associated primes.
	
	The following lemma is quite useful, and it shows how studying absorbing ideals will help us better understand the ideals of Noetherian rings.
	\begin{lemma}\textup{\cite[Corollary 6.8]{n-abs}.}\label{Lemma: omega alpha is finite for Noetherian rings}
		For an ideal $\alpha$ of a Noetherian ring $A$, \[\omega(\alpha)\leq\omega^*(\alpha)<\infty\]
		\begin{proof}
			Given that $\omega(\alpha)\leq\omega^*(\alpha)$ is always true, it is enough to show $\omega^*(\alpha)<\infty$.
			We will prove this by following Anderson and Badawi's proof \cite[Corollary 6.8]{n-abs}. 
			Since $A$ is Noetherian, $\alpha$ must have a primary decomposition $\alpha=\bigcap_{i=1}^r Q_i$. 
			We know $\omega^*(\bigcap_{i=1}^r Q_i)\leq \sum_{i=1}^r \omega^*(Q_i)$ and hence it suffices to show $\omega^*(Q)<\infty$ for any $\p$-primary ideal $Q$ of $A$. 
			Since $Q$ is $\p$-primary, we must have $\omega^*(Q)= \omega_{A_\p}^*(Q_\p)$. 
			Passing to the quotient ring gives us $\omega_{A_\p}^*( Q_\p)=\omega_{A_\p/ Q_\p}^*(0)$.
			Notice that $A_\p/ Q_\p := A_\p/A_\p Q$ is a zero dimensional Noetherian ring, and thus a local Artin ring.
			Hence there exists some finite $k$ such that $(\p_\p/Q_\p)^k=0$, and so $\omega_{A_\p/ Q_\p}^*(0)\leq k$. 
			This implies $\omega^*(Q)<\infty$, as needed.
		\end{proof}
	\end{lemma}
	
	We will now present a simplified version of Choi and Walker's solution for \hyperlink{(C3)}{(C3)}\textemdash modifying, concatenating, and shortening arguments for brevity. 
	It is recommend that the reader see \cite{Choi_Walker_arxiv} or \cite{Choi_Walker} for insight into how this proof was obtained.
	\begin{lemma}\textup{\cite[Theorem 1]{Choi_Walker}.}\label{Choi-Walker}
		If $\omega(\alpha)<\infty$ then \[(\sqrt{\alpha})^{\omega(\alpha)}\seq \alpha\]
		\begin{proof}
			Let $n=\omega(\alpha)$. Notice that $(\sqrt{\alpha})^n\seq \alpha$ if and only if $J^n\seq \alpha$ for every subideal $J\seq\sqrt{\alpha}$ generated by $n$ elements $a_1,...,a_n\in \sqrt{\alpha}$. 
			For any $a_i\in \sqrt{\alpha}$ we know $a_i^n\in \alpha$, and thus $J^{n^2}=\<\{a_i\}_{i}\>^{n^2}\seq \<\{a_i^n\}_i\>\seq\alpha$. 
			Let $\mathcal{N}$ be the set of vectors $v\in \mathbb{N}^n$ such that $n\leq\sum_i v_i \leq n^2$. Notice that $J^n \seq \alpha$ if and only if $\prod_ia_i^{v_i}\in\alpha$ for every $v$ in $\mathcal{N}$. 
			The finite set $\mathcal{N}$ is ordered lexicographically, and so we can use reverse induction to show $\prod_ia_i^{v_i}\in\alpha$ for all $v\in \mathcal{N}$. 
			The base case $({n^2},0,...,0)$ is trivial since $a_1^{n^2}\in J^{n^2}\seq \alpha$. 
			Fix $r\in \mathcal{N}$ and assume $\prod_ia_i^{v_i}\in\alpha$ for every $v\in \mathcal{N}$ lexicographically greater than $r$.
			
			Let $g=\prod_i a_i^{r_i}$, we will show $g\in\alpha$. Let $g_{j} :=a_j^{r_j-1}\prod_{i\ne j} a_i^{r_i}$ for each $j$ with $r_j\ne 0$.  
			If $\sum_i r_i=n^2$ then $g\in J^{n^2}\seq \alpha$ and were done, assume not.
			By the induction hypothesis we have $ga_i\in \alpha$ for all $a_i$, and thus  $g y\in \alpha$ for every $y\in J=\<\{a_i\}_i\>$. 
			Yet $g=\prod_i a_i^{r_i}$ and $\omega(\alpha)=n\leq\sum_i r_i$, so $gy\in \alpha$ implies either $g\in \alpha$ or $g_j y\in\alpha$ for some $r_j\ne 0$. 
			If $g\in\alpha$ we are done, so assume that for every $y\in J$ there exists $r_j\ne 0$ with $g_j y\in\alpha$. 
			
			Let $j_0$ be maximal such that $r_{j_0}\ne 0$, and define the sequence $\{j_i\}_{i=0}^\infty$ recursively by letting $j_{k}$ be largest such that $r_{j_k}\ne 0$ and $g_{j_{k}}\sum_{i<k} a_{j_i}\in\alpha$. 
			If $j_k<j_{k+1}$ for some $k$, then by the induction hypothesis $g_{j_{k+1}} a_{j_k}\in \alpha$. Hence
			\[g_{j_{k+1}} \sum_{i<k} a_{j_i} = g_{j_{k+1}} \sum_{i<k+1} a_{j_i}  -g_{j_{k+1}}a_{j_k}\in\alpha\]
			and thus $g_{j_{k+1}} \sum_{i<k} a_{j_i} \in\alpha$, a contradiction by choice of $j_k$ largest. 
			Therefore $j_{i+1}\leq j_i$ for all $i$.
			Hence $\{j_i\}_i$ is both infinite and non-increasing with $j_i\in\{1,...,k\}$, and thus $j_{k+1}= j_k$ for some $k$. Since $g_{j_{k}} \sum_{i<k} a_{j_i}\in\alpha$ and $g = g_{j_{k}}a_{j_k}$, we know $g+\alpha=g_{j_{k}}a_{j_k}+g_{j_{k}} \sum_{i<k} a_{j_i}+\alpha$. This implies that:
			\[	g+\alpha= g_{j_{k}}a_{j_k}+g_{j_{k}} \sum_{i<k} a_{j_i}+ \alpha = g_{j_{k}} \sum_{i<k+1} a_{j_i}+ \alpha = g_{j_{k+1}} \sum_{i<k+1} a_{j_i}+ \alpha =\alpha	\]
			which shows $g\in \alpha$. Thus $J^n\seq\alpha$, and therefore $(\sqrt{\alpha})^n\seq \alpha$.
		\end{proof}
	\end{lemma}
	
	Before we give a new proof that \hyperlink{(C2)}{(C2)} implies \hyperlink{(C1)}{(C1)}, we will need some definitions and a lemma. 
	The content $c(f)$ of a polynomial $f(x)=\sum_i a_i x^i\in A[x]$ is the ideal $c(f)=\<\{a_i\}_i\>$ of the ring $A$ generated by the coefficients of $f$. 
	Given an ideal $\beta$ of $A[x]$, the content $c(\beta)$ of $\beta$ is the ideal given by $c(\beta)=\sum_{f\in \beta} c(f)$, which is the smallest ideal of $A$ such that $\beta\seq c(\beta)[x]$. In particular, $\beta\seq\alpha[x]$ if and only if $c(\beta)\seq\alpha$.
	
	From these definitions we get the following trivial lemma, whose proof is left as an exercise for the reader.
	\begin{lemma}\label{equation: key fact for splitting content of polynomials}
		For any two polynomials $g,h\in A[x]$, if $\deg(g)<\ell$ then:
		\[c(g)\cdot c(h)=c(g(x)\cdot h(x^\ell))\quad \text{and}\quad c(g)+ c(h)=c(g(x)+x^\ell h(x)) \]
	\end{lemma}
	
	Donadze \cite[Proposition 2.7]{Donadze_infinite_fields} was the first to publish a proof that \hyperlink{(C2)}{(C2)} implies \hyperlink{(C1)}{(C1)}, although Laradji independently proved the same result \cite{Laradji}.
	Remarkably, Laradji submitted his first draft exactly one day before Donadze did\textemdash but their proofs differ greatly.
	Laradji's approach used the following stronger result, which we will prove using a different technique.
	
	\begin{claim}\label{omega* alpha is leq omega alpha[x]}\textup{\cite[Proposition 2.10]{Laradji}.} We always have
		\(\omega^*(\alpha)\leq \omega(\alpha[x])\).
		\begin{proof} 
			Suppose $\omega(\alpha[x])=n$, and recall that $\omega^*(\alpha)=n$ if and only if $\alpha$ is strongly $n$-absorbing for finitely generated ideals. Let $[\gamma_i]_{i=1}^{n+1}$ be a list of finitely generated ideals $\gamma_i=\<\{a_{ij}\}_{j=0}^{n_i}\>$ with $\prod_i \gamma_i\seq \alpha$. 
			We will show that there exists $t$ such that $\prod_{i\ne t} \gamma_i \seq \alpha$.
			
			Let $k$ be larger than the size of each of the finite sets $\{a_{ij}\}_j$, i.e. $k>\max_i(n_i+1)$.
			For each $i$, let $f_i = \sum_j a_{ij} x^{k^i\cdot j}\in A[x^{k^{i}}]$ and notice that $c(f_i)=\gamma_i$. 
			Hence 
			\(c\left(\prod_i f_i \right)\seq \prod_i c(f_i) =\prod_i\gamma_i \seq \alpha\).
			We know $\prod_i f_i \in \alpha[x]$ if and only if $c\left(\prod_i f_i \right)\seq \alpha$, and thus
			\(\prod_i f_i \in \alpha[x]\).
			Since $\alpha[x]$ is n-absorbing, there must exist some $t$ such that 
			\( \prod_{i\ne t} f_i \in \alpha[x] \).
			
			By choice of $k$, for each $i$ we know $k$ is larger than the size of $\{a_{ij}\}_j$ and thus $\deg(f_i)<k\cdot k^{i}=k^{i+1}$. Yet $f_j\in A[x^{k^{i+1}}]$ for each $j>i$, and so
			by repeatedly applying Lemma (\ref{equation: key fact for splitting content of polynomials}) we get $c(\prod_{i\ne t}f_i)=\prod_{i\ne t}c(f_i)$. Since $c(f_i)=\gamma_i$ for all $i$, we have $\prod_{i\ne t}\gamma_i=c(\prod_{i\ne t}f_i)\seq \alpha$.
		\end{proof}
	\end{claim}
	
	Before we move forward, let us discuss an example of a non-Noetherian ring whose ideals are all $2$-absorbing. 
	This example shows that non-Noetherian rings can satisfy the statement of Lemma (\ref{Lemma: omega alpha is finite for Noetherian rings}).
	\begin{example}\label{Example: the F(+)V construction}
		The technique we give here is inspired by that seen in \cite[Example 3.4]{n-abs} and \cite[Theorem 3.4]{2-abs}.
		For a $R$-module $M$, the \textit{idealization} of $M$ is the ring $R(+)M$ defined as the $R$-module $R\oplus M$ with a multiplicative operation given by $(r,m)(s,n)=(rs,rn+sm)$. 
		The origin of this construction is unclear, although the first systematic study appears in Nagata's book \cite{nagata1962local}.
		This book coined the term \enquote{principle of idealization} in reference to the method of studying a module $M$ by considering the ideals of $R(+)M$.
		
		Let $F$ be a field, let $V$ be a vector space over $F$, and consider the ring $A=F(+) V$.
		For any $(u,m)\in A$ with $u\ne 0$, we have $(u,m)(u^{-1},-u^{-2}m)=(1,0)$ and thus the non-units of $A$ are of the form $(0,m)$; in particular, the subspaces $W$ of $V$ are in a natural one-to-one correspondence with the ideals $0(+)W$ of $A$.
		This implies that $A$ is Noetherian if and only if $V$ is finite dimensional, in which case $A$ is also Artinian.
		Choosing $V$ infinite dimensional guarantees non-Noetherianness.
		
		We will show that every ideal $\alpha$ of $A$ is $2$-absorbing.
		Suppose $\prod_i (r_i,m_i)\in \alpha$. If $r_i=r_j=0$ for some $i\ne j$, then $(r_i,m_i)(r_j,m_j)=(0,m_i)(0,m_j)=(0,0)\in \alpha$ and were done. 
		Assume not.
		Since non-units are of the form $(0,m)$, and since $\prod_i (r_i,m_i)=(\prod_i r_i,\sum_j (\prod_{i\ne j} r_i) m_j )\in \alpha$, we know $\prod_i r_i=0$. 
		Therefore, there has to be exactly one element satisfying $r_k=0$, which implies $(r_i,m_i)$ is a unit for any $i\ne k$. 
		Hence $(r_k,m_k)\in \alpha$ and so $\omega(\alpha)\leq 2$.
		Therefore, every ideal of $A$ will be $2$-absorbing.
		Since fields are the only rings whose ideals are all prime, this is the strongest such example we can give.
	\end{example}

	\section{Collections of an Ideal}\label{section: collections of an ideal}
	Before we continue, let us first clarify and review some notation.
	For two ideals $\alpha$ and $\beta$ of a ring $A$, the colon ideal $(\alpha:\beta)$ is the ideal defined by
	\((\alpha:\beta)=\{f\in A\;|\;f\beta\seq\alpha\}\).
	In other words, $(\alpha:\beta)$ is the largest ideal of $A$ such that $(\alpha:\beta)\cdot \beta\seq \alpha$. 
	We will always have $(\alpha:\beta\gamma)=((\alpha:\beta):\gamma)$, and for any ideals $\alpha\seq\alpha'$ and $\beta\seq \beta'$ we have $(\alpha:\beta')\seq(\alpha':\beta)$.
	Moreover, if $S$ is a multiplicatively closed subset, then $S^{-1}(\alpha:\beta)\seq (S^{-1}\alpha:S^{-1}\beta)$ with equality whenever $A$ is Noetherian or $\beta$ is finitely generated.
	By a slight abuse of notation, we let $(\alpha:f)=(\alpha:\<f\>)$ when $f\in A$, and $(\alpha:W)=(\alpha:\<W\>)$ when $W$ is a subset of $A$. 

	Given finite lists $L=[L_1,...,L_n]$ and $L'=[L_1',...,L_k']$, we will use \[L\cup L' = [L_1,...,L_n,L'_1,...,L'_k]\] to represent the concatenation of the list $L$ by $L'$. e.g. $[3]\cup [2,3]=[3,2,3]$.
	
	For a finite list $C=[C_1,...,C_n]$ we will use $|C|$ to represent the size of $C$, in this case $|C|=n$. We write $C_{i\ne j}$ to represent the sublist of $C$ which misses the $j^{th}$ component. i.e. \[C_{i\ne j}=[C_1,...,C_{j-1},C_{j+1},...,C_n]\]
	Similarly, given $I=\{i_1,...,i_k\}\seq \{1,...,n\}$ we will write $C_{i\in I}$ for the list $[C_{i_1},...,C_{i_k}]$, and  write $C_{i\not\in I}$ to represent $C_{i\in \{1,...,n\}\setminus I}$.

	\subsection{Maximal and Semimaximal Collections}\label{subsection: ontro to maximal collections}
	Recall every maximal collection is semimaximal, and a collection $N$ is semimaximal if and only if $N_j=(\alpha:\prod_{i\ne j}N_i)$ for all $j$.
	So if $N$ is a maximal or semimaximal collection of $\alpha$, then $\alpha\seq N_j$ for all $N_j\in N$. 
	Hence for any subideal $\beta\seq\alpha$, the maximal collections of $\alpha/\beta$ are in one-to-one correspondence with the maximal collections of $\alpha$, and the same is true for semimaximal collections. 
	
	By definition of $\omega^*(\alpha)$, any semimaximal $\omega^*(\alpha)$-collection of $\alpha$ must also be a maximal collection. 
	This is not true in general; there are semimaximal collections which are not maximal.
	For instance, consider the semimaximal $2$-collection $[\<2\>,$ $\<4\>]$ of the ideal $\<8\>$ in the ring $\Z$. This collection is not maximal, as $[\<2\>,$ $\<4\>]$ lies below the (unique) maximal collection $[\<2\>,$ $\<2\>,$ $\<2\>]$ of $\<8\>$. 
	Moreover, $[\alpha]$ is always a semimaximal $1$-collection of $\alpha$, but it is only a maximal collection when $\alpha$ is a prime ideal.
	
	\begin{examples}
		Let $\alpha$ be an ideal of a Dedekind domain $\mathcal{O}$, and let $\alpha=\prod_{i=1}^n \p_i$ be the factorization of $\alpha$ into primes, repeating copies of the same prime ideal when needed. Up to reordering, $[\p_i]_i$ is the unique maximal collection of $\alpha$, and the semimaximal collections are of the form $\left[\prod_{j\in I_i}\p_j\right]_i$ for some \textit{partition} $\bigcup_i I_i=\{1,...,n\}$.
		
		Now let $\m=0(+)V$ be the maximal ideal of the the ring $A=F(+)V$, as discussed in Example (\ref{Example: the F(+)V construction}). 
		Since $\m$ is prime, we know the maximal collection $[\m]$ is the only semimaximal collection of $\m$. 
		We have $\m^2=\<0\>$ and thus $[\m,\m]$ will be the unique maximal collection of any ideal $\alpha\ne \m$.
		Moreover, the only other semimaximal collection of $\alpha$ will be the trivial collection $[\alpha]$.
		We will need more machinery before we can discuss less trivial examples in depth. 
	\end{examples}

	Maximal and semimaximal collections have a lot of structure which will prove useful later on. 
	Before we discuss this, let us prove a result about the existence of these collections.
	\begin{lemma}\label{Lemma: Maximal collection existence} 
		Let $\alpha$ be an ideal of a ring $A$. Every $n$-collection of $\alpha$ lies below at least one semimaximal $n$-collection. If $A$ is Noetherian, or if $\omega^*(\alpha)<\infty$, then every collection of $\alpha$ lies below some maximal collection of $\alpha$.
		\begin{proof}
			By Zorn's Lemma, every $n$-collection must be contained in some semimaximal $n$-collection.
			
			By Lemma (\ref{Lemma: omega alpha is finite for Noetherian rings}) we know $\omega^*(\alpha)<\infty$ for any ideal $\alpha$ of a Noetherian ring. Thus we need only prove that if $\omega^*(\alpha)<\infty$ then every collection lies below a maximal collection. Let $C$ be a collection of $\alpha$. There must exist some $k\leq \omega^*(\alpha)$ such that $C$ lies below a collection $C'$ of size $k$ but $C$ does not lie below any collection of size greater than $k$. By the previous part of this proof, $C'$ must lie below a semimaximal $k$-collection $N$. Since \enquote{lying below} is a transitive relation, and since $C$ cannot lie below a collection of size greater than $k$, we know $N$ must be a maximal collection. 
			Since $C\leq N$, we can conclude that $N$ is a maximal collection which lies above $C$.	
		\end{proof}
	\end{lemma}
	If $\omega^*(\alpha)=\infty$ then maximal collections may not even exist. 
	Indeed, consider the ideal $\alpha=\<x^\pi\>$ of the monoid ring $\Z[x^{\mathbb{R}^{\geq 0}}]\cong\Z[(\mathbb{R}^{\geq 0},+)]$. 
	Notice that $\prod_{i=1}^n x^{\pi/n}\seq \alpha$ defines a collection of $\alpha$ for any $n$, and thus $\omega^*(\alpha)=\infty$.
	The semimaximal collections of $\alpha$ will be of the form $N=[\<x^{a_i}\>]_i$ where $\sum_i a_i=\pi$. 
	These collections will never be maximal since $N$ will always lie below the collection $N_{i\ne j}\cup[\<x^{{a_j}/{2}}\>,\<x^{{a_j}/{2}}\>]$ for any $N_j=\<x^{a_j}\>$.
	Since all maximal collections are semimaximal, we know $\alpha$ cannot have any maximal collections.
	
	In preparation for the next lemma, let us review some facts about associated primes. For an ideal $\alpha$ of a Noetherian ring $A$, the set of \textit{associated primes} $\ass(\alpha)$ is the (finite) set of prime ideals $\p$ which satisfy
	$\p=(\alpha:f)$ for some $f\in A$. 
	It is well know that for any primary decomposition $\alpha=\bigcap_i Q_i$ we have $\ass(\alpha)=\{\sqrt{Q_i}\}_i$.
	
	Let us show that when $A$ is Noetherian, $\p$ is an associated prime of $\alpha$ if and only if $\p=(\alpha:(\alpha:\p))$. 
	The \enquote{only if} direction is trivial, and if $\p=(\alpha:(\alpha:\p))$ then since $A$ is Noetherian we know $\<\{f_i\}_i\>= (\alpha:\p)$ for some finite set $\{f_i\}_i$. 
	We have $\p=(\alpha:(\alpha:\p))=(\alpha:\<\{f_i\}_i\>)=\bigcap_i (\alpha:f_i)$, and since $\{f_i\}_i$ is finite we know $\p=(\alpha:f_i)$ for some $f_i$.
	Many of the results connecting collections and associated primes can be generalized to the non-Noetherian context by considering prime ideals satisfying $\p=(\alpha:(\alpha:\p))$, although we will not need this generality here.
	
	The following lemma shows that any associated prime of any element of a semimaximal collection of $\alpha$ must also be an associated prime of $\alpha$. 
	This fact will become useful when computing specific examples (see the discussion preceding (\ref{example: classifying maximal collections})).
	
	\begin{lemma}\label{Lemma: associated primes of a maximal collection are a subset of ass alpha}
		Let $\alpha$ be an ideal of a Noetherian ring, and let $N$ be a semimaximal collection of $\alpha$. Then $\bigcup_i \ass(N_i) \seq \ass(\alpha)$, where $\ass(\alpha)$ is the set of associated primes of $\alpha$.
		\begin{proof}
			Since $N_j=(\alpha:\prod_{i\ne j}N_i)$ for any element $N_j$ of a semimaximal collection $N$, it suffices to show that $\ass(\alpha:\gamma)\seq \ass(\alpha)$ for any ideal $\gamma\not\seq \alpha$.
			By the discussion preceding this lemma, $\p$ is an associated prime of $(\alpha:\gamma)$ if and only if $\p=((\alpha:\gamma):((\alpha:\gamma):\p))$. Assume $\p$ is an associated prime of $(\alpha:\gamma)$.
			Since $((\beta_1:\beta_2):\beta_3)=(\beta_1:\beta_2\beta_3)$ for any ideals $\beta_i$, we get
			\[\p=((\alpha:\gamma):((\alpha:\gamma):\p))=(\alpha:\gamma(\alpha:\gamma \p))\]
			and thus $\p=(\alpha:\gamma(\alpha:\gamma \p))$. Clearly $\gamma(\alpha:\gamma \p)\seq(\alpha:\p)$, and hence 
			\[(\alpha:(\alpha:\p))\seq (\alpha:\gamma(\alpha:\gamma \p))=\p\]
			This implies $(\alpha:(\alpha:\p))\seq\p$, but clearly $\p\seq (\alpha:(\alpha:\p))$ and thus $\p=(\alpha:(\alpha:\p))$. Therefore, any associated prime of $(\alpha:\gamma)$ will also be an associated prime of $\alpha$. 
		\end{proof}
	\end{lemma}
	
	Now let us focus on maximal collections. 
	As we will see in (\ref{example: classifying maximal collections}), the ideal $\<x,y,z\>^2\cap\<x^2,y\>$ of $\mathbb{Q}[x,y,z]$ has the maximal collection $[\<x^2,y\>,$ $\<x,y,z\>]$, which contains the non-prime element $\<x^2,y\>$. 
	Therefore, not all elements of maximal collections are prime ideals. 
	In spite of this, we will prove a slightly weaker statement: that every \textit{maximal} element of a maximal collection is a prime ideal.
	
	\begin{lemma}\label{Lemma: Max elements of max collections are prime} 
		Let $\alpha$ be an ideal of some ring $A$.
		\begin{enumerate}
			\item The maximal elements of a maximal collection are all prime ideals. 
			\item If $M_k$ is an element of a maximal collection $M$, then for any non-trivial collection $C$ of $M_k$, there must exist a maximal element $M_m$ of $M$ such that $\prod_{i\ne j}C_i\seq M_m$ for all $j$, and $\prod_iC_i\prod_{i\ne k,m}M_i\seq \alpha$. 
			\item Suppose $A$ is Noetherian, let $M$ be a maximal collection of $\alpha$, and let $\p$ be a maximal element of the finite set of prime ideals $\bigcup_i \ass(M_i)$. Then $\p\in \ass(M_i)$ if and only if $\p=M_i$. 
		\end{enumerate}	
		\begin{proof}
			Clearly (i) follows from (ii), so we will first prove (ii). Suppose $M=[M_1,...,M_n]$ is a maximal collection of $\alpha$, and let $M_k$ be an element of $M$. Notice that $M_k$ is a prime ideal, if and only if, $[M_k]$ is the only maximal collection of $M_k$. Suppose $C=[C_1,...,C_\ell]$ is a collection of $M_k$, and without loss of generality we can let $C$ be semimaximal. This implies $M_k\seq C_i$ for any $C_i\in C$, and thus
			\[M\leq M_{i\ne k} \cup C\]  
			By maximality of $M$, we know $M_{i\ne k} \cup C$ cannot be a collection of $\alpha$. 
			Therefore, either $\prod_{i\ne m} C_i\prod_{i\ne k} M_i\seq\alpha$ for some $m$, or $\prod_{i}C_i\prod_{i\ne k,j} M_i\seq\alpha$ for some $j\ne k$. 
			If $\prod_{i\ne m}C_i\prod_{i\ne k} M_i\seq\alpha$ then since $M$ is maximal it must be semimaximal, and thus we have \(\prod_{i\ne m}C_i\seq (\alpha:\prod_{i\ne k}M_i)=M_k\). This is a contradiction since $C$ is a collection of $M_k$. Hence $\prod_{i}C_i\prod_{i\ne k,j} M_i\seq\alpha$ for some $j\ne k$. 
			
			First we will show that we can choose $M_j$ such that $M_j$ is a maximal element of $M$.
			If $M_j\seq M_k$ then since $M_k\seq C_1$ we would have $\prod_{i\ne 1}C_i\prod_{i\ne k,j} M_i\seq\prod_{i}C_i\prod_{i\ne k,j} M_i\seq\alpha$, a contradiction since $\prod_{i\ne m} C_i\prod_{i\ne k} M_i\not\seq\alpha$ for all $m$.
			Now if $M_j\seq M_{j'}$ for some $j'\ne k$, then $\prod_{i}C_i\prod_{i\ne k,j'} M_i\seq\prod_{i}C_i\prod_{i\ne k,j} M_i\seq\alpha$. Therefore, we can choose $j$ such that $M_j$ is a maximal element of $M$.		
			For each $m$ we have $M_k\seq C_m$, and so
			\[\prod_{i\ne m}C_i\cdot\prod_{i\ne j}M_i= \prod_{i\ne m}C_i \cdot M_k\cdot \prod_{i\ne j,k}M_i\seq \prod_{i\ne m}C_i \cdot C_m\cdot \prod_{i\ne j,k}M_i\seq \alpha\]
			Hence $\prod_{i\ne m}C_i\seq (\alpha: \prod_{i\ne j}M_i)$ for all $m$. Since $M$ is maximal, we know $M_j=(\alpha: \prod_{i\ne j}M_i)$ which implies $\prod_{i\ne m}C_i\seq M_j$ for all $m$. This finishes the proof of (ii).
			
			Finally, we will prove (iii). Suppose $M_k\in M$ such that $\p\in \ass(M_k)$ for some maximal element $\p$ of $\bigcup_i \ass(M_i)$. With intent to contradict, assume that $M_k\ne \p$. Then $(M_k:\p)\ne M_k$, and thus $[\p,(M_k:\p)]$ must be a 2-collection of $M_k$. By (ii) there exists some $m$ such that both $\p$ and $(M_k:\p)$ are contained in $M_m$. 
			Yet $\p $ is a maximal element of $\bigcup_i \ass(M_i)$, and clearly $M_m\seq \q$ for all $\q\in\ass(M_m)$. Since $\p\seq M_m$ and by maximality of $\p$, we know $\ass(M_m)=\{\p\}$, and thus $M_m=\p$. Yet we know that $\p(M_k:\p)\prod_{i\ne k,m}M_i\seq \alpha$, and hence $M_m(M_k:\p)\prod_{i\ne k,m}M_i\seq \alpha$. 
			This implies that $(M_k:\p)\seq(\alpha:\prod_{i\ne k}M_i)$. 
			Since $M$ is maximal, we know $(\alpha:\prod_{i\ne k}M_i)=M_k$ and thus $(M_k:\p)\seq M_k$. This implies that $M_k=(M_k:\p)$, and thus $\p\not\in \ass(M_k)$\textemdash a contradiction. 
			Therefore $M_k=\p$, as needed.
		\end{proof}
	\end{lemma}
	
	Now we will briefly shift focus to generating examples. 
	Consider a Noetherian ring $A$, and let $\alpha$ be an ideal of $A$ with primary decomposition $\alpha=\bigcap_{i=1}^k Q_i$. When the size of $\ass(\alpha)$ and values of $\omega(Q_i)$ are both very small, it becomes possible to use lemmas (\ref{Lemma: Max elements of max collections are prime}) and (\ref{Lemma: associated primes of a maximal collection are a subset of ass alpha}) to classify the maximal collections of $\alpha$. Indeed, every maximal collection $M$ of $\alpha$ must contain an element $M_j$ with $M_j\in\ass(\alpha)$, and thus $M_{i\ne j}$ will be a semimaximal collection of $(\alpha:M_j)$. Therefore, if we can classify all semimaximal collections of $(\alpha:\p)$ for every $\p\in \ass(\alpha)$, then we can find every maximal collection of $\alpha$. For example:
	
	\begin{examples}\label{example: classifying maximal collections}
		Consider the ring $\mathbb{Q}[x,y,z]$ and the ideal $\alpha=\<x,y,z\>^2\cap\<x^2,y\>$. We will show that, up to reordering, $\alpha$ has two maximal collections: $[\<x^2,y\>,$ $\<x,y,z\>]$ and $[\<x,y\>,$ $\<x,y\>]$. 
		By Lemma (\ref{Lemma: Max elements of max collections are prime}), in order to find all the maximal collections of $\alpha$, it is enough to find all semimaximal collections of $(\alpha:\p)$ for every associated prime $\p$ of $\alpha$. 
		The only two associated primes of $\alpha$ are $\<x,y\>$ and $\<x,y,z\>$. Clearly, $(\alpha:\<x,y,z\>)=\<x^2,y\>$ has two semimaximal collections up to reordering: $[\<x^2,y\>]$ and $[\<x,y\>,$ $\<x,y\>]$. The list $[\<x,y\>,$ $\<x,y\>,$ $\<x,y,z\>]$ is not a collection since $\<x,y\>^2\seq\alpha$, so we need only consider $[\<x^2,y\>,$ $\<x,y,z\>]$, which is clearly semimaximal since $(\alpha:\<x^2,y\>)=\<x,y,z\>$ and $(\alpha:\<x,y,z\>)=\<x^2,y\>$.
		For the element $\<x,y\>\in\ass(\alpha)$, notice that $(\alpha:\<x,y\>)=\<x,y\>$ is a prime ideal and thus will only have one semimaximal collection $[\<x,y\>]$, which gives us a semimaximal collection $[\<x,y\>,$ $\<x,y\>]$ of $\alpha$.
		Therefore reorderings of $[\<x,y\>,$ $\<x,y\>]$ and $[\<x^2,y\>,$ $\<x,y,z\>]$ are the only semimaximal collections of $\alpha$ which contain a prime ideal. 
		Neither lies above the other, so up to reordering $\alpha$ has two maximal collections: $[\<x^2,y\>,$ $\<x,y,z\>]$ and $[\<x,y\>,$ $\<x,y\>]$.
		
		Now let $\beta = \<x,y,z\>^2\cap \<x^2,y\>\cap \<a,x^3,y,z\>$ be an ideal of the polynomial ring $\mathbb{Q}[a,x,y,z]$. By following the same process, we get that the only semimaximal collections of $\beta$ which contain a prime ideal are reorderings of the two collections $[\<x,y\>,$ $\<x,y\>,$ $\<a,x,y,z\>]$ and $[\<x^2,y\>,$ $\<x,y,z\>]$. Since $[\<x^2,y\>,$ $\<x,y,z\>]$ lies below the collection $[\<x,y\>,$ $\<x,y\>,$ $\<a,x,y,z\>]$, we know that the only maximal collections of $\beta$ are reorderings of $[\<x,y\>,$ $\<x,y\>,$ $\<a,x,y,z\>]$.
	\end{examples}
	
	Now we will discuss localizations of semimaximal collections.
	Given a multiplicatively closed subset $S$ of $A$, and a collection $C$ of an ideal $\alpha$ of $A$, define the localization of $C$ by $S$ to be the list $S^{-1}C=[S^{-1}C_i]_i$.
	Notice that $S^{-1}C$ is often not a collection. 
	For example, $[\<2\>,$ $\<4\>,$ $\<3\>]$ is a semimaximal collection of the ideal $\<24\>$ of $\Z$, yet the localization $(\Z\setminus  \<2\>)^{-1} [\<2\>,$ $\<4\>,$ $\<3\>] = [\<2\>_{\<2\>},$ $\<4\>_{\<2\>},$ $\<1\>_{\<2\>}]$ contains the unit ideal. 
	Nonetheless, the nonunit elements $[\<2\>_{\<2\>}, $ $\<4\>_{\<2\>}]$ form a semimaximal collection of $\<24\>\Z_{\<2\>}=\<8\>_{\<2\>}$.
	This is true for any semimaximal collection of a Noetherian ring, as seen by the following lemma.
	
	\begin{lemma}\label{Lemma: max collection correspondences}
		Let $\alpha$ be an ideal of a Noetherian ring $A$, and let $S$ be a multiplicatively closed subset of $A$ disjoint from $\alpha$. Then we have:
		\begin{enumerate}
			\item If $M$ is a maximal collection of $\alpha$ such that $S\cap \bigcup_i M_i=\emptyset$, then $S^{-1}M$ is a maximal collection of $S^{-1}\alpha$.
			\item For any semimaximal collection $N$ of $\alpha$, the sublist containing all the nonunit elements of $S^{-1}N$ will be a semimaximal collection of the ideal $S^{-1}\alpha$. 
			Conversely, for every semimaximal collection $N$ of $S^{-1}\alpha$ there exists a semimaximal collection $N'$ of $\alpha$ such that $N$ will be equal to the sublist of $S^{-1}N'$ containing the non-unit elements of $S^{-1}N'$.
			\item When $n=\omega^*(\alpha)<\infty$ the maximal $n$-collections of $S^{-1}\alpha$ will correspond to maximal $n$-collections $M$ of $\alpha$ with $S\cap \bigcup_i M_i = \emptyset$.
		\end{enumerate}
		\begin{proof}
			Let us prove (i). 
			Clearly $S^{-1}M$ will be a collection of $S^{-1}$.
			Suppose $M'$ is a maximal collection of $S^{-1}\alpha$ lying above $S^{-1}M$, and let $\pi:A\rightarrow S^{-1}A$ be canonical.
			We will first show $\prod_i\pi^{-1}(M'_i)\seq \alpha$. Indeed, if $\prod_i\pi^{-1}(M'_i)\not\seq \alpha$ then by Noetherianness of $A$ we have $S\cap(\alpha:\prod_i\pi^{-1}(M'_i))\ne \emptyset$.
			It is easy to check that $[\pi^{-1}(M'_i)]_i\cup[(\alpha:\prod_i\pi^{-1}(M'_i))]$ forms a collection of $\alpha$ which lies above $M$, and thus
			by maximality $M=[\pi^{-1}(M'_i)]_i\cup[(\alpha:\prod_i\pi^{-1}(M'_i))]$.
			This is a contradiction since $\emptyset=(\alpha:\prod_i\pi^{-1}(M'_i))\cap S\seq \bigcup_i M_i\cap S = \emptyset$.
			Hence $\prod_i\pi^{-1}(M'_i)\seq \alpha$ and thus $[\pi^{-1}(M'_i)]_i$ forms a collection which lies above $M$. 
			By maximality of $M$ we have $M=[\pi^{-1}(M'_i)]_i$, and thus $M'=S^{-1}M$.
			
			Now we will prove (ii). 
			Suppose $N$ is a semimaximal collection of $\alpha$, and let $N'$ be the sublist of $N$ containing the elements $N_i$ with $S\cap N_i=\emptyset$.
			Since $A$ is Noetherian, for any $N_{j'}=N_j'\in N'$ we have 
			\[\<1\>\ne S^{-1}N_j=S^{-1}(\alpha:\prod_{i\ne j'}N_i)=(S^{-1}\alpha:S^{-1}\prod_{i\ne j'}N_i)=(S^{-1}\alpha:\prod_{i\ne j}S^{-1}N_i')\]
			and therefore $\prod_{i\ne j}S^{-1}N_i'\not\seq S^{-1}\alpha$ and $S^{-1}N_j'=(S^{-1}\alpha:\prod_{i\ne j}S^{-1}N_i')$. 
			It follows that $S^{-1}N'$ must forms a semimaximal collection of $S^{-1}\alpha$; in particular, the nonunit elements of $S^{-1}N$ form a semimaximal collection of $S^{-1}\alpha$.
			
			Suppose $N'$ is a semimaximal collection of $S^{-1}\alpha$, and let $N_i=\pi^{-1}(N'_i)$. 
			Notice that $S^{-1}N_i=N'_i$ and $(\alpha:\prod_{i\ne j}N_i)\seq\pi^{-1}(S^{-1}\alpha:\prod_{i\ne j}N'_i)=\pi^{-1}(N'_i)=N_i$. 
			Thus if $\prod_i N_i\seq \alpha$ then $N=[N_i]_i$ is a collection of $\alpha$, and since $(\alpha:\prod_{i\ne j}N_i)\seq N_i$ we know $N$ would be a semimaximal collection. 
			Hence we may assume that $\prod_i N_i\not\seq \alpha$. 
			
			Notice that $(\alpha:\prod_i N_i)\ne \<1\>$ but $S\cap (\alpha:\prod_i N_i)\ne \emptyset$. 
			Thus $S^{-1}((\alpha:\prod_i N_i) \cdot \prod_{i\ne j}N_i)=\prod_{i\ne j}N_i'\not\seq S^{-1}\alpha$, and so $[(\alpha:\prod_i N_i)]\cup [N_i]_i$ is a collection of $\alpha$. 
			Since $N'$ is semimaximal, $N'_j=(S^{-1}\alpha:\prod_{i\ne j}N'_i)$ for all $j$, and so
			\[\pi^{-1}(S^{-1}\alpha:S^{-1}(\alpha:\prod_i N_i)\prod_{i\ne j}N_i)=\pi^{-1}(S^{-1}\alpha:\prod_{i\ne j}N'_i)=\pi^{-1}(N'_j)=N_j\]
			which implies $(\alpha:(\alpha:\prod_i N_i)\prod_{i\ne j}N_i)\seq N_j$ and therefore $N_j=(\alpha:(\alpha:\prod_i N_i)\prod_{i\ne j}N_i)$ for all $j$. 
			We can now conclude that $[(\alpha:\prod_i N_i)]\cup [N_i]_i$ is a semimaximal collection of $\alpha$, as needed.
			
			(iii) follows from (ii), and thus we are done.
		\end{proof}
	\end{lemma}
	
	Unfortunately, if $S\cap \bigcup_i M_i\ne\emptyset$ then the nonunit elements of $S^{-1}M$ need not form a maximal collection of $S^{-1}\alpha$. Moreover, not every maximal collection of $S^{-1}\alpha$ can be equal to the set of nonunit elements of $S^{-1}M$ for any maximal collection $M$ of $\alpha$.
	The following examples demonstrate this.
	
	\begin{examples}\label{example: localizing maximal collections}
		Consider the ideal $\alpha=\<x,y,z\>^2\cap\<x^2,y\>$ of $\mathbb{Q}[x,y,z]$. By (\ref{example: classifying maximal collections}), the maximal collections of $\alpha$ are reorderings of either $[\<x^2,y\>,$ $\<x,y,z\>]$ or $[\<x,y\>,$ $\<x,y\>]$. 
		The localization of $[\<x,y\>,$ $\<x,y\>]$ by the prime ideal $\<x,y\>$ is the the maximal collection $[\<x,y\>_{\<x,y\>},$ $\<x,y\>_{\<x,y\>}]$. We know this localization must be a maximal collection because every element of $[\<x,y\>,$ $\<x,y\>]$ is disjoint with $S:=A\setminus \<x,y\>$.
		In contrast, the localization of $[\<x^2,y\>,$ $\<x,y,z\>]$ is the list $[\<x^2,y\>_{\<x,y\>},$ $\<1\>]$, 
		and clearly $[\<x^2,y\>_{\<x,y\>}]$ is a semimaximal collection which lies below the (unique) maximal collection $[\<x,y\>_{\<x,y\>},$ $\<x,y\>_{\<x,y\>}]$ of $\alpha_{\<x,y\>}=\<x^2,y\>_{\<x,y\>}$. 
		Therefore, the nonunit elements of the localization of a maximal collection need not form a maximal collection in general.
		
		Now consider larger polynomial ring $\mathbb{Q}[a,x,y,z]$, and let $\beta = \<x,y,z\>^2\cap \<x^2,y\>\cap \<a,x^3,y,z\>$. 
		The localization of $\beta$ by the prime ideal $\<x,y,z\>$ will have the maximal collection $[\<x^2,y\>_{\<x,y,z\>},$ $\<x,y,z\>_{\<x,y,z\>}]$. 
		By our discussion in (\ref{example: classifying maximal collections}), we know the only maximal collections of $\beta$ are reorderings of $[\<x,y\>,$ $\<x,y\>,$ $\<x,y,z,a\>]$, and so the maximal collection $[\<x^2,y\>_{\<x,y,z\>},$ $\<x,y,z\>_{\<x,y,z\>}]$ cannot equal the nonunit elements of $(A\setminus\<x,y,z\>)^{-1}M$ for any maximal collection $M$ of $\beta$.
		Therefore, not every maximal collection of $S^{-1}\beta$ can be equal to the nonunit elements of $S^{-1}M$ for some maximal collection $M$ of $\beta$.	
	\end{examples} 
	
	If $\omega^*(\alpha)=n$ then every semimaximal $n$-collection of $\alpha$ must also be a maximal collection. Thus, as a corollary of Lemma (\ref{Lemma: max collection correspondences}), we have:
	
	\begin{lemma}\label{Lemma: ideal extension from smallest Ui Mi must exist and be n absorbing}
		Let $\alpha$ be an ideal of a Noetherian ring $A$ with $\omega^*(\alpha)=n$. Then there exists at least one maximal $n$-collection $M$ such that $\bigcup_i M_i$ is minimal amongst all maximal $n$-collections of $\alpha$.
		
		Moreover, if we set $S:=A\setminus \bigcup_i M_i$ then $\omega^*(S^{-1}\alpha)=n$ and every maximal $n$-collection of $S^{-1}\alpha$ will contain a copy of each maximal ideal of $S^{-1}A$.
		\begin{proof}
			Since $A$ is Noetherian, we know $\ass(\alpha)$ will be a finite set. By lemmas (\ref{Lemma: associated primes of a maximal collection are a subset of ass alpha}) and (\ref{Lemma: Max elements of max collections are prime}), we know that the maximal elements of any maximal collection of $\alpha$ must be a subset of the finite set $\ass(\alpha)$. So there must exist some maximal $n$-collection $M$ such that $\bigcup_i M_i$ is minimal amongst maximal $n$-collections. 
			
			The maximal elements of $M$ are all prime ideals, and thus $S:=A\setminus\bigcup_i M_i$ must be a multiplicatively closed subset of $A$. 
			If $\p_1,...,\p_k$ are the maximal elements of $M$, then $S=A\setminus \bigcup_{i=1}^k \p_i$ and so by the Prime Avoidance Lemma the maximal ideals of $S^{-1}A$ will be of the form $S^{-1}\p_i$. 	
			
			Since all semimaximal $n$-collections are maximal $n$-collections, we can use Lemma (\ref{Lemma: max collection correspondences}) to conclude that the maximal $n$-collections $S^{-1}M'$ of $S^{-1}\alpha$ correspond to the maximal $n$-collections $M'$ of $\alpha$ whose maximal elements are all contained in $\bigcup_i M_i=\bigcup_{i=1}^k \p_i$. 
			By choice of $M$, for any maximal $n$-collection $S^{-1}M'$ of $S^{-1}A$ we have $\bigcup_iS^{-1}M'_i=\bigcup_iS^{-1}M_i=\bigcup_{i=1}^k S^{-1}\p_i$.
			By Lemma (\ref{Lemma: Max elements of max collections are prime}), the maximal elements of $S^{-1}M'$ must be prime ideals. Using the Prime Avoidance Lemma, we know $S^{-1}M'$ must contain at least one copy of each maximal element of $S^{-1}M$, which are precisely the maximal ideals of $S^{-1}A$.
			This finishes the proof.
		\end{proof}
	\end{lemma}

	Before we move on to the next subject, we should briefly discuss an interesting correspondence --- although we will not use it here. 
	Let $T(A)$ be the total ring of fractions, which is the localization of $A$ by the set of normal elements $S=\{f\in A\;|\;\<0\>=(\<0\>:f)\}$. 
	By definition of $S$, we know $S^{-1}(\<0\>:\beta)=(S^{-1}\<0\>:S^{-1}\beta)$.
	Moreover, if $\pi:A\rightarrow T(A)=S^{-1}A$ is canonical then $\pi^{-1}((S^{-1}\<0\>:\beta'))=(\<0\>:\pi^{-1}(\beta'))$ for any ideal $\beta'$ of $T(A)$. 
	It follows that the semimaximal (resp. maximal) collections of the zero ideal in $A$ correspond naturally to the semimaximal (resp. maximal) collections of the zero ideal in $T(A)$.
	In particular, the maximal collections of any ideal $\alpha$ are in a one-to-one correspondence with the maximal collections of the zero ideal in the total ring of fractions $T(A/\alpha)$, and the same is true of semimaximal collections.

	\subsection{The Tilde Ideals $\tilde C_I$}\label{subsection: tilde ideals}
	In this subsection, we introduce the tilde ideals $\tilde C_I$ of a collection $C$. Studying these ideals can give us insight into the structure of a collection $C$, and tilde ideals provide a convenient setting to generalize some of the results we have on semimaximal collections.
	
	Let $\alpha$ be an ideal of some ring $A$, and let $C=[C_1,...,C_k]$ be a collection of $\alpha$. 
	Given a subset $I\seq \{1,...,k\}$ we define the \textit{tilde ideal} $\tilde C_I$ of $C$ to be:
	\[\tilde C_I := (\alpha:\prod_{i\not\in I} C_i) = \{x\in A\;|\;x\prod_{i\not\in I} C_i\seq \alpha\}\]
	We have $\tilde C_\emptyset =\<1\>$, and for $I\ne \emptyset$ we know $\tilde C_I$ is the largest ideal such that $C_{i\not\in I}\cup[\tilde C_I]$ is a collection of $\alpha$. 
	By abuse of notation, we write $\tilde C_j=\tilde C_{\{j\}}$ and $\tilde C_{j_1,...,j_m}=\tilde C_{\{j_1,...,j_m\}}$.
	Clearly $C_j\seq \tilde C_j$, and furthermore, the list $C_{i\in I}$ is always a collection of $\tilde C_I$.
	We know $N$ is a semimaximal collection if and only if $N_j=\tilde N_j:=(\alpha:\prod_{i\ne j}N_i)$ for all $j$. 
	In particular, if $|N|=\omega^*(\alpha)$ then $N$ is a maximal collection if and only if $N_j=\tilde N_j$ for all $N_j$.
	
	For a subideal $\beta$ of $\alpha$, let $C/\beta = [C_1/\beta,...,C_k/\beta]$ be the induced collection of $\alpha/\beta$. We always have:
	\[\tilde C_I/\beta = \tilde{(C/\beta)}_I\]
	For a multiplicatively closed set $S$, the list $S^{-1}C:=[S^{-1}C_i]_i$ does not need to be a collection of $S^{-1}\alpha$. 
	As a remedy for this, the following generalization of Lemma (\ref{Lemma: max collection correspondences}) shows that tilde ideals can be used to describe the sublists $S^{-1}C_{i\not\in I}:=[S^{-1}C_i]_{i\not\in I}$ of $S^{-1}C$ which are collections of $S^{-1}\alpha$, whenever $\alpha$ is an ideal of a Noetherian ring.

	\begin{claim}\label{prop: localization and tilde ideals}
		For an ideal $\alpha$ of a Noetherian ring $A$ and a multiplicatively closed subset $S$ disjoint from $\alpha$, the sublists $S^{-1}C_{i\not\in I}\seq S^{-1}C$ which are collections of $S^{-1}\alpha$ correspond to the maximal subsets $I\seq\{1,...,k\}$ with $S\cap\tilde C_{I}\ne \emptyset$.
		\begin{proof}
			First, let us generalize the tilde ideal construction to finite lists $L$ containing ideals $L_i$ such that $\prod_iL_i\seq \alpha$. We can do this by simply setting $\tilde L_I=(\alpha:\prod_{i\not\in I}L_i)$. 
			By construction of $\tilde L_{I}$, the sublists $L_{i\not\in I}:=[L_i]_{i\not\in I}$ of $L$ which are collections of $\alpha$ correspond naturally to the maximal subsets $I\seq\{1,...,k\}$ satisfying $\tilde L_{I}= \<1\>$.
			Notice that the list $S^{-1}C:=[S^{-1}C_i]_i$ has the property that $\prod_i S^{-1}C_i\seq S^{-1}\alpha$. 
			Since $A$ is Noetherian, we have $S^{-1}\tilde C_I = (S^{-1}\alpha:\prod_{i\not\in I}S^{-1}C_i)=\tilde{(S^{-1}C)}_I$ and thus $\tilde{(S^{-1}C)}_I=\<1\>$ if and only if $S\cap\tilde C_I\ne\emptyset$. 
			The claim follows.
		\end{proof}
	\end{claim}
	
	For an ideal $\alpha$ of a Noetherian ring, the following generalization of Lemma (\ref{Lemma: associated primes of a maximal collection are a subset of ass alpha}) shows that $\ass(\tilde C_I)\seq \ass(\alpha)$ for any collection $C$ and set of indices $I$.
	
	\begin{lemma}\label{Lemma: associated primes of tilde ideals}
		Let $\alpha$ be an ideal of a Noetherian ring $A$. For any collection $C$ of $\alpha$, we have $\bigcup_{I} \ass(\tilde C_I)\seq \ass(\alpha)$.
		\begin{proof}
			This follows from the proof of Lemma (\ref{Lemma: associated primes of a maximal collection are a subset of ass alpha}).
		\end{proof}
	\end{lemma}
	
	With the above lemma in mind, the following is a weak generalization of (iii) from Lemma (\ref{Lemma: Max elements of max collections are prime}) for collections of largest possible length.
	
	\begin{lemma}\label{Lemma: ideal prod has a max ideal as an associated prime then it is equal to that maximal ideal}
		Let $\alpha$ be a strongly $n$-absorbing ideal of a Noetherian ring $A$, and let $C$ be an $n$-collection of $\alpha$. If $\ass(\tilde C_j)$ contains a maximal element $\p$ of $\ass(\alpha)$, then $\tilde C_j=\p$.
		\begin{proof}
			Notice that $[\tilde C_j]\cup C_{i\ne j}$ is an $n$-collection of $\alpha$. Let $M$ be a maximal collection of $\alpha$ lying above $[\tilde C_j]\cup C_{i\ne j}$. Since $\alpha$ is strongly $n$-absorbing, $M$ must be an $n$-collection. 
			Without loss of generality, we can reorder $M$ such that $\tilde C_j\seq M_j$ and $C_i\seq M_i$ for all $i\ne j$. Hence $\prod_{i\ne j} C_i\seq \prod_{i\ne j} M_i$ and thus 
			\(\tilde M_j=(\alpha:\prod_{i\ne j} M_i)\seq (\alpha:\prod_{i\ne j} C_i)=\tilde C_j\).
			Yet $M$ is maximal and thus semimaximal, so $M_j=(\alpha:\prod_{i\ne j} M_i)=\tilde M_j$ and thus $M_j=\tilde C_j$.
			Since $\p$ is a maximal element of $\ass(\alpha)$ and $\p\in \ass( M_j)$, by Lemma (\ref{Lemma: associated primes of a maximal collection are a subset of ass alpha}) we know $\p$ must be a maximal element of $\bigcup_i\ass(M_i)$.
			By (iii) of Lemma (\ref{Lemma: Max elements of max collections are prime}), we have $M_j=\p$ and therefore $\tilde C_j=\p$.
		\end{proof}
	\end{lemma}

	\subsection{Maximal Collections in Artin Rings}\label{subsection: artin rings}
	If $\alpha$ is an ideal of a ring $A$, then the maximal collections of $\alpha$ will correspond canonically to the maximal collections of the zero ideal in the ring $A/\alpha$. 
	In particular, if $A/\alpha$ is Artinian then we can use facts about the maximal collections of the zero ideal in an Artin ring to describe the maximal collections of $\alpha$. 
	As we will see in Lemma (\ref{base case}), such an ideal $\alpha$ will have a unique maximal collection up to reordering.
	
	It is a well known fact that any Artin ring $A$ is isomorphic to the (finite) product of its localizations $A\cong \prod_\m A_\m$. Hence the following will be useful when proving Lemma (\ref{base case}).
	
	\begin{lemma}\label{Lemma: product rings}
		Let $A=\prod_{i=1}^r A_i$ be a product ring. Then,
		\begin{enumerate}
			\item $M$ is a semimaximal collection of $\<0\>_{A}$, if and only if, for each $k$ there exists some semimaximal collection $M^{(k)}=[M_j^{(k)}]_j$ of $\<0\>_{A_k}$ such that (after reordering) we have
			\[M = \left[\prod_{i< k}A_i\times M^{(k)}_j \times \prod_{i> k}A_i\right]_{k,j} = \bigcup_k \left[\prod_{i< k}A_i\times M^{(k)}_j \times \prod_{i> k}A_i\right]_{j}\]
			Moreover, $M$ is a maximal collection of $\<0\>_{A}$ if and only if each $M^{(k)}$ is a maximal collection of $\<0\>_{A_k}$.
			\item $\omega^*_{\prod_i A_i}(0)=\sum_i \omega^*_{A_i}(0)$ and $\omega_{\prod_i A_i}(0)=\sum_i \omega_{A_i}(0)$.
		\end{enumerate} 
		\begin{proof}
			(i) is trivial, and (ii) follows from (i).
		\end{proof}
	\end{lemma}
	
	Now let us prove the following lemma, which we will use in the proof of our main result: Theorem (\ref{Theorem: (C1) - the First Anderson-Badawi Conjecture}).
	
	\begin{lemma}\label{base case}
		Let $A$ be an Artinian ring. Up to reordering, there exists a unique maximal collection $M$ of $\<0\>$. Moreover, each $M_i\in M$ must be a maximal ideal of $A$, and for each $i$ there exists $f_i\in M_i$ such that $\prod_i f_i=0$ and $\prod_{i\ne j}f_i\ne0$ for all $j$.
		\begin{proof} 
			First we will prove this for local Artinian rings $A_\m$. Indeed, since any Artin ring is Noetherian, by Lemma (\ref{Lemma: omega alpha is finite for Noetherian rings}) we know that $\omega_{A_\m}(0)<\infty$. By Proposition (\ref{Choi-Walker}) we know $\m^{\omega_{A_\m}(0)}=\<0\>$. Since $A_\m$ is local, we know that $\m^{\omega_{A_\m}(0)}=\<0\>$ implies $\omega_{A_\m}^*(0)\leq \omega_{A_\m}(0)$. Yet $\omega_{A_\m}(0)\leq \omega_{A_\m}^*(0)$ and thus $\omega_{A_\m}(0)=\omega_{A_\m}^*(0)$. It follows trivially that the only maximal collection of $A_\m$ is the collection consisting of $\omega_{A_\m}(0)$ copies of $\m$.
			
			Now that we have proven this lemma for local Artinian rings, lets generalize this to arbitrary Artinian rings. 
			Every Artinian ring $A$ is canonically isomorphic to the product of its localizations $A\cong\prod_i A_{\m_i}$, and thus by Lemma (\ref{Lemma: product rings}) we know that every maximal collection of $A$ will be a reordering of the collection $M$ consisting of $\omega_{A_\m}(0)$ copies of $\m$ for each maximal ideal $\m$. Again by Lemma (\ref{Lemma: product rings}), we have $\omega_{A}(0)=\sum_\m \omega_{A_\m}(0)=\sum_\m \omega_{A_\m}^*(0)=\omega^*_A(\alpha)$. 
			In particular, there must exist a principal $\omega^*(\alpha)$-collection $[\<f_i\>]_i$. Since $[\<f_i\>]_i$ must lie below some maximal collection of $A$, and any maximal collection of $A$ is a reordering of $M$, we can let $f_i\in M_i$ for all $i$, which finishes the proof.
		\end{proof}
	\end{lemma}

	\section{$n$-Absorbing Ideals are Strongly $n$-Absorbing}\label{section: the first conjecture}
	In this section we will prove our main result, Theorem (\ref{Theorem: (C1) - the First Anderson-Badawi Conjecture}). This Theorem states that $\omega(\alpha)=\omega^*(\alpha)$ for any ideal $\alpha$, as conjectured \hyperlink{(C1)}{(C1)} in \cite{n-abs}. 
	Before we prove this, we will need two well-known results; one from Neal McCoy and one from Edward Davis.
	
	The following argument is a fundamental result in the study of finite unions of ideals (see \cite{Gottlieb}), which first appears as an unnamed lemma inside the proof of the first theorem in McCoy's infamous paper \cite{McCoy}. 
	This lemma will serve a critical role in Lemma (\ref{Lemma: reduction to base case}).
	
	\begin{lemma}\label{Lemma: McCoy's ideal union lemma}\textup{\cite{McCoy}.}
		Let $\alpha, \gamma_1,...,\gamma_n$ be ideals of a ring $A$ such that $\alpha\seq \bigcup_i \gamma_i$ but $\alpha\not\seq \bigcup_{i\ne j}\gamma_i$ for all $j$. Then for any $j$ we have
		\[\alpha\cap \bigcap_{i\ne j} \gamma_i\seq \gamma_j\]
		\begin{proof}
			We will follow McCoy \cite[Proof of Theorem 1]{McCoy}.
			Since $\alpha\not\seq \bigcup_{i\ne j}\gamma_i$ we know there must exist some $y\in \gamma_j\cap \alpha\setminus \bigcup_{i\ne j}\gamma_i$. 
			Given any $x\in \alpha\cap\bigcap_{i\ne j}\gamma_i$ we must have $x+y\in \alpha$, and hence there must be some $i$ such that $x+y\in \gamma_i$. 
			For any $i\ne j$, we know $x\in \gamma_i$ and $y\not\in \gamma_i$, and therefore $x+y\not\in \gamma_i$. 
			This implies $x+y\in \gamma_j$, yet $y\in \gamma_j$ and thus $x\in \gamma_j$. By arbitrary choice of $x$, we know $\alpha\cap \bigcap_{i\ne j}\gamma_i\seq \gamma_j$.
		\end{proof}
	\end{lemma}
	
	Recall that the Prime Avoidance Lemma states that for any finite set of primes $\{\p_i\}_i$ and any ideal $\alpha$, if $\alpha\seq\bigcup_i \p_i$ then $\alpha\seq\p_i$ for some $i$.
	The following generalization of the Prime Avoidance Lemma appears in \cite[Exercise 16.8]{matsumura1989commutative}, where Matsumura credits the result to Edward Davis (without citation).
	This result is often referred to as \enquote{Davis' Prime Avoidance Lemma} due to its surprising utility, as we will see in the proof of Theorem (\ref{Theorem: (C1) - the First Anderson-Badawi Conjecture}).
	
	\begin{lemma}\textup{\textbf{(Davis' Prime Avoidance Lemma).}}\label{Lemma: Davis' Prime Avoidance}
		Suppose $\beta$ is an ideal of $A$, let $\{\p_i\}_i$ be a finite set of prime ideals of $A$, and suppose $f\in A$. If $\<f\>+\beta\not\seq\p_i$ for all $i$, then there exists some $g\in \beta$ such that $f+g\not\in \bigcup_i \p_i$.
		\begin{proof}
			Assume $\p_j\not\seq \p_i$ for all $j\ne i$, and let $I$ be the set with $i\in I$ whenever $f\in \p_i$. 
			We must have $\beta\not\seq \p_i$ for all $i\in I$, since if not then $\<f\>+\beta\not\seq \p_i$ implies $f\not\in \p_i$ and so $i\not\in I$.
			Thus $\beta\cdot \prod_{j\not\in I}\p_j\not \seq \p_i$ for all $i\in I$, and 
			by the Prime Avoidance Lemma, $\beta\cdot \prod_{j\not\in I}\p_j\not\seq\bigcup_{i\in I} \p_i$.
			Choosing any $g\in \beta\cdot \prod_{j\not\in I}\p_j\setminus\bigcup_{i\in I} \p_i$ gives us an element $g\in \beta$ with $f+g\not\in \bigcup_i\p_i$.
		\end{proof}
	\end{lemma}

	\subsection{The Main Theorem}\label{subsection: n-absorbing ideals are strongly n-absorbing}
	In this subsection, we will use our previous lemmas to complete the final steps in our proof that $\omega(\alpha)=\omega^*(\alpha)$.
	The following lemma is a key step in this proof, as it allows us to reduce the problem to a convenient setting.
	
	\begin{lemma}\label{Lemma: reduction to quasilocal noetherian rings}
		$\omega(\alpha)=\omega^*(\alpha)$ for any ideal $\alpha$, if and only if, $\omega(\alpha)=\omega^*(\alpha)$ for any ideal $\alpha$ of a quasilocal Noetherian ring such that every maximal $\omega^*(\alpha)$-collection of $\alpha$ contains a copy of each maximal ideal of the ring.
		\begin{proof}
			The \enquote{only if} direction is tautological. We will prove the other direction using a contrapositive argument. Assume $\alpha$ is an ideal of a ring $A$ such that $\omega(\alpha)<\omega^*(\alpha)$. We will show that there exists a finitely generated subring $B$ of $A$ and a multiplicatively closed subset $S\seq B$ such that $\omega(S^{-1}(\alpha\cap B))<\omega^*(S^{-1}(\alpha\cap B))$ and the ring $S^{-1}B$ is a quasilocal Noetherian ring with the property that every maximal $\omega^*(S^{-1}(\alpha\cap B))$-collection of $S^{-1}(\alpha\cap B)$ contains a copy of each maximal ideal of $S^{-1}B$.
			
			We know $\omega(\alpha)<\omega^*(\alpha)$, so there must exist a collection $C$ of $\alpha$ with $|C|$ greater than $\omega(\alpha)$. 
			We have $\prod_{i\ne j} C_i\not\seq \alpha$ for each $j$, and hence for each $i\ne j$ there must exist some $f_{ij}\in C_i$ such that $\prod_{i \ne j} f_{ij}\not\in \alpha$. 
			It follows that $\<\{f_{ij}\}_{j\ne i}\>\seq C_i$ and thus the list $[\<\{f_{ij}\}_{j\ne i}\>]_i$ is a collection of size greater than $\omega(\alpha)$, where $\{f_{ij}\}_{j\ne i}$ iterates over $j$ for a fixed $i$. 
			Let $B$ be the finitely generated subring of $A$ generated by each $f_{ij}$. 
			Similarly, $[\,\<\{f_{ij}\}_{j\ne i}\>_B\,]_i$ is a collection of $\alpha\cap B$ of size greater than $\omega(\alpha)$. Since $\omega(\alpha\cap B)\leq \omega(\alpha)$, we know $\omega(\alpha\cap B)<\omega^*(\alpha\cap B)$.
			Therefore, we have a finitely generated ring $B$ with $\omega(\alpha\cap B)<\omega^*(\alpha\cap B)$.
			
			Since $B$ is a finitely generated ring, $B$ must be Noetherian, and so by Lemma (\ref{Lemma: omega alpha is finite for Noetherian rings}) we know $\omega^*(\alpha\cap B)<\infty$. 
			Let $n=\omega^*(\alpha\cap B)$.
			By Lemma (\ref{Lemma: ideal extension from smallest Ui Mi must exist and be n absorbing}) there must exist a maximal $n$-collection $M$ of $\alpha\cap B$ such that $\bigcup_i M_i$ is minimal. 
			Let $S=B\setminus\bigcup_i M_i$ and since Lemma (\ref{Lemma: Max elements of max collections are prime}) shows that the maximal elements of $M$ are prime ideals of $B$, we know $S$ is a multiplicatively closed subset. By the Prime Avoidance Theorem, $S^{-1}B$ must be a quasilocal Noetherian ring. 
			Again by Lemma (\ref{Lemma: ideal extension from smallest Ui Mi must exist and be n absorbing}), we know $\omega^*(S^{-1}(\alpha\cap B))=\omega^*(\alpha\cap B)$ and any maximal $n$-collection of $S^{-1}(\alpha\cap B)$ will contain a copy of each and every maximal ideal of $S^{-1}B$. 
			Since $\omega(S^{-1}(\alpha\cap B))\leq \omega(\alpha\cap B)< \omega^*(\alpha\cap B)=\omega^*(S^{-1}(\alpha\cap B))$, we know $\omega(S^{-1}(\alpha\cap B))<\omega^*(S^{-1}(\alpha\cap B))$ as needed.
		\end{proof}
	\end{lemma}
	
	With Lemma (\ref{Lemma: reduction to quasilocal noetherian rings}) in mind, let us prove the final and most important lemma of Theorem (\ref{Theorem: (C1) - the First Anderson-Badawi Conjecture}).
	
	\begin{lemma}\label{Lemma: reduction to base case}
		Let $A$ be a quasilocal Noetherian ring, and let $\alpha$ be an ideal of $A$ with $\omega^*(\alpha)=n$ such that every semimaximal $n$-collection of $\alpha$ contains each and every maximal ideal of $A$. Then there must exist some $n$-collection $C$ such that for every $C_j\in C$ either 
		\begin{enumerate}
			\item $\tilde C_j$ is a maximal ideal and $C_j=\tilde C_j$, or,
			\item $\tilde C_j$ is not a maximal ideal and $C_j$ is a principal ideal.
		\end{enumerate} 
		And moreover, $C$ contains a copy of every maximal ideal of $A$.
		\begin{proof}
			Consider an $n$-collection $C$ of $\alpha$ such that $C$ contains the largest possible number of elements $C_j$ which are maximal ideals of $A$. 
			Since $[\tilde C_j]\cup C_{i\ne j}$ is also an $n$-collection, we know whenever $\tilde C_j$ is a maximal ideal we must have $C_j=\tilde C_j$.
			Now if $C_k$ is not a maximal ideal and $f\in \tilde C_k$ such that $[\<f\>]\cup C_{i\ne k}$ is a collection of $A$, then $C'=[\<f\>]\cup C_{i\ne k}$ must contain every element of $C$ which is a maximal ideal. 
			If $\tilde C'_j$ is a maximal ideal, then $C'_j=\tilde C'_j$ since if not then $[\tilde C'_j]\cup C'_{i\ne j}$ is a $n$-collection of $\alpha$ which contains more copies of maximal ideals than $C$ does. 
			Therefore, if we can show that whenever $C_k$ is not maximal or principal, there must exist $f\in C_k$ such that $[\<f\>]\cup C_{i\ne k}$ is a collection of $A$, then we can repeat this process to get an $n$-collection whose elements satisfy either (i) or (ii).
			
			Before we continue, we must show that $C$ contains a copy of every maximal ideal of $A$. 
			Indeed, by Lemma (\ref{Lemma: Maximal collection existence}) there must be a maximal collection $M$ lying above $C$ with $C_i\seq M_i$ for all $i$. Since $|C|=\omega^*(\alpha)$, we know $|M|=|C|$, and therefore $M_i=\tilde M_i \seq \tilde C_i$ for all $i$. 
			We have assumed that every maximal $\omega^*(\alpha)$-collection must contain a copy $M_{i}=\m$ of every maximal ideal $\m$. 
			Thus we have $\m=M_{i} \seq \tilde C_{i}$ and hence $\tilde C_{i}=\m$.
			By choice of $C$, we know $C_{i}=\tilde C_{i}=\m$, and therefore $C$ must contain a copy of each maximal ideal $\m$ of $A$.
			In particular, the maximal elements of $C$ are maximal ideals.
			
			Suppose there is some $C_k\in C$ which is both non-maximal and non-principal. 
			Recall from $\S$\ref{subsection: tilde ideals} that $\tilde C_{k,j}:=(\alpha:\prod_{i\ne k,j}C_i)$.
			If $C_k\not\seq\bigcup_{j\ne k} \tilde C_{k,j}$ then there exists some $f\in C_k$ such that $f\not\in  \bigcup_{j \ne k}\tilde C_{k,j}= \bigcup_{j \ne k}(\alpha:\prod_{i\ne j,k}C_i)$, and hence $\<f\>\prod_{i\ne k,j}C_i\not \seq \alpha$ for any $j\ne k$. 
			Since $f\in C_k$, we know $\<f\>\prod_{i\ne k}C_i \seq \alpha$. 
			Yet $C$ is a collection, so $\prod_{i\ne k}C_i\not \seq \alpha$. 
			It follows that $[\<f\>]\cup C_{i\ne k}$ must be a collection of $\alpha$. 
			Thus we need only prove that $C_k\not\seq\bigcup_{i\ne k} \tilde C_{k,i}$.
			
			With intent to contradict, assume that $ C_k\seq \bigcup_{i\ne k}\tilde C_{k,i}$. 
			We will show this implies there exists $\ell\ne k$ such that $ C_k\seq  \tilde C_{k,\ell}$, which is a contradiction since $\tilde C_{k,\ell}:=(\alpha:\prod_{i\ne k,\ell}C_i)$ and thus we would have $\prod_{i\ne \ell}C_i\seq \alpha$.
			
			Whenever $C_j\seq C_t$ for some $j$ and $t$ not equal to $k$, then we have $\prod_{i\ne k,t}C_i\seq \prod_{i\ne k,j}C_i$ and thus 
			$ \tilde C_{k,j} = (\alpha : \prod_{i\ne k,j}C_i)\seq (\alpha : \prod_{i\ne k,t}C_i)=\tilde C_{k,t}$. 
			Since every maximal element of $C$ is a maximal ideal, and since $C_k$ is not maximal, there must exist some positive integers $m_1,...,m_r\ne k$ such that each $C_{m_i}$ is a maximal ideal of $A$ with $ C_k\seq \bigcup_{i=1}^r\tilde C_{k,m_i}$ and $ C_k\not\seq \bigcup_{i\ne j}\tilde C_{k,m_i}$ for all $j$. 
			It follows that $\tilde C_{k,m_i}\ne \tilde C_{k,m_j}$ for all $i\ne j$, and thus $C_{m_i}\ne  C_{m_j} $ for all $i\ne j$. 
			
			By Lemma (\ref{Lemma: McCoy's ideal union lemma}) we must have $C_k\cap \bigcap_{i\ne j}\tilde C_{k,m_i} \seq \tilde C_{k,m_j}$ for any $j$.
			Yet by definition of $\tilde C_{k,m_i}$ we have $ C_k C_{m_i}\seq\tilde C_{k,m_i}$ for any $i$. 
			Hence for any $j$ we have
			\[C_k \cdot \bigcap_{{i\ne j}} C_{m_i} \seq  \bigcap_{i\ne j} C_k C_{m_i}\seq C_k\cap \bigcap_{i\ne j}\tilde C_{k,m_i}\seq \tilde C_{k,m_j}\]
			This implies $C_k \cdot \bigcap_{{i\ne j}} C_{m_i}\seq \tilde C_{k,m_j}$ and therefore $\bigcap_{{i\ne j}} C_{m_i}\seq (\tilde C_{k,m_j}:C_k)$.
			Since $(\tilde C_{k,m_j}:C_k)=\tilde C_{m_j}$ and $\tilde C_{m_j}=C_{m_j}$, we know $\bigcap_{{i\ne j}} C_{m_i}\seq C_{m_j}$.
			This is a contradiction since $C_{m_j}$ is a prime ideal and $C_{m_i}\not\seq C_{m_j}$ for all $i\ne j$, which implies $\bigcap_{{i\ne j}} C_{m_i}\not\seq C_{m_j}$. 
			Hence $ C_k\not\seq \bigcup_{i\ne k}\tilde C_{k,i}$, as needed. 
		\end{proof}
	\end{lemma}
	
	Finally, we can prove our main result.
	We will do this by first reducing to the setting of Lemma (\ref{Lemma: reduction to quasilocal noetherian rings}), where we can create a principal $\omega^*(\alpha)$-collection by combining lemmas (\ref{Lemma: reduction to base case}) and (\ref{base case}) using Davis' Prime Avoidance Lemma.
	
	\begin{theorem}\label{Theorem: (C1) - the First Anderson-Badawi Conjecture}
		Let $\alpha$ be any ideal of an arbitrary ring $A$. Then, \[\omega(\alpha)=\omega^*(\alpha)\]
		\begin{proof}
			By Lemma (\ref{Lemma: reduction to quasilocal noetherian rings}), we may assume that $A$ is a quasilocal Noetherian ring and $\alpha$ is an ideal of $A$ with $\omega^*(\alpha)=n$ such that every semimaximal $n$-collection of $\alpha$ contains each and every maximal ideal of $A$. 
			Now let $C$ be an $n$-collection of $\alpha$ as in Lemma (\ref{Lemma: reduction to base case}), let $[M_i]_i$ be the sublist of $C$ containing all maximal elements of $C$, and let $[\<h_i\>]_i$ be the sublist of $C$ containing all non-maximal elements of $C$. 
			
			After reordering we have $C=[M_i]_i\cup [\<h_i\>]_i$, and thus $[M_i]_i$ must be a collection of the ideal $\beta:=(\alpha:\prod_i h_i)$; i.e. $\beta=\tilde C_I$ where $I$ is such that $C_{i\in I}=[M_i]_i$. 
			Since the elements of $[M_i]_i$ are all maximal ideals of $A$, we know that $[M_i]_i$ must be a maximal collection of $\beta$. 
			Since $\prod_i M_i\seq \beta\seq M_i$, we know that the maximal ideals of $A/\beta$ must be elements of the set $\{M_i/\beta\}_i$.
			
			Since $\prod_i (M_i/\beta)=0/\beta$, we know every prime ideal of $A/\beta$ must contain some maximal ideal $M_i/\beta$, and thus $A/\beta$ must be Artinian.
			By Lemma (\ref{base case}), we know that for each $i$ there must exist some $f_i\in M_i$ such that $[\<f_i\>]_i$ is a collection of $\beta$. 
			Since $\prod_i M_i\seq \beta$, for any $\p\in \ass(\alpha)\setminus\{ M_i\}_i$ we know $\beta\not\seq \p$ and thus $\<f_j\>+\beta\not\seq \p$ for any $j$. 
			By Davis' Prime Avoidance Lemma (\ref{Lemma: Davis' Prime Avoidance}), for each $f_j$ there is some $g_j\in \beta$ such that $f_j+g_j\not\in \bigcup (\ass(\alpha)\setminus \{M_i\}_i)$.
			Hence without loss of generality, we may assume $f_j\not\in \bigcup (\ass(\alpha)\setminus \{M_i\}_i)$ for all $j$.
			
			We will show that $\prod_i h_i \prod_k f_k\in \alpha$ defines a principal $n$-collection. Indeed, since $[\<f_k\>]_k$ forms a collection of $\beta=(\alpha:\prod_i h_i)$ we know that $\prod_{k\ne j}f_k \not\seq (\alpha:\prod_i h_i)$, and thus $\prod_i h_i\prod_{k\ne j}f_k \not\seq \alpha$ for all $j$. Hence we need only show that $\prod_{i\ne j} h_i\prod_kf_k \not\seq \alpha$ for all $j$.
			Let $S=A\setminus \bigcup (\ass(\alpha)\setminus \{M_i\}_i)$, and notice that $S$ is multiplicatively closed. 
			Since $f_j\not\in \bigcup (\ass(\alpha)\setminus \{M_i\}_i)$ for all $j$, we know $f_j\in S$ for all $j$. 
			Moreover, $M_j\cap S\ne \emptyset$ for all $j$. This implies:
			\[	(S^{-1}\alpha:S^{-1}\<\prod_{i\ne j} h_i\prod_kf_k\>)
			= (S^{-1}\alpha:S^{-1}\<\prod_{i\ne j} h_i\>) 
			= (S^{-1}\alpha:S^{-1}\prod_{i\ne j} h_i\prod_k M_k) \]
			and therefore $S^{-1}(\alpha:\prod_{i\ne j} h_i\prod_kf_k) = S^{-1}(\alpha:\prod_{i\ne j} h_i\prod_k M_k) $.
			Let $j'$ be such that $C_{j'}=\<h_j\>$, and thus $\tilde C_{j'}=(\alpha:\prod_{i\ne j} h_i\prod_k M_k)$.
			By Lemma (\ref{Lemma: associated primes of tilde ideals}), we know $\ass(\tilde C_{j'})=\ass(\alpha:\prod_{i\ne j} h_i\prod_k M_k)\seq \ass(\alpha)$. 
			By Lemma (\ref{Lemma: ideal prod has a max ideal as an associated prime then it is equal to that maximal ideal}) and choice of $C$,
			we know $\ass(\tilde C_{j'})\seq \ass(\alpha)\setminus \{M_i\}_i$. Thus 
			\[\<1\>\ne  S^{-1}\tilde C_{j'} =S^{-1}(\alpha:\prod_{i\ne j} h_i\prod_k M_k) =S^{-1}(\alpha:\prod_{i\ne j} h_i\prod_kf_k) \]
			which implies $(\alpha:\prod_{i\ne j} h_i\prod_kf_k) \ne \<1\>$ and therefore $\prod_{i\ne j} h_i\prod_kf_k\not\in \alpha$ for every $j$. 
			Hence $\prod_i h_i \prod_k f_k\in \alpha$ defines a principal $n$-collection $[\<h_i\>]_i\cup[\<f_k\>]_k$ of $\alpha$, as needed.
		\end{proof}
	\end{theorem}

	\section{Some New Conjectures}\label{section: some new conjectures}
	In this section, we introduce and discuss three new conjectures that are closely related to \hyperlink{(C2)}{(C2)}\textemdash the last unsolved conjecture on absorbing ideals given by Anderson and Badawi. 
	The relationship between these conjectures can be summarized in the following diagram:
	\begin{center}
		\begin{tikzcd}[remember picture, overlay]
			\begin{array}{c}
				\textup{\text{(\ref{conjecture: finitely many maximal collections})}}
			\end{array}
			\arrow[rrr,Rightarrow, "\text{Prop. } (\ref{prop: conjecture finitely many maximal collections implies conjecture maximal collections of alpha[x] are M[x]})
			"']
			&[0em] &[-2em] &[0em]
			\begin{array}{c}
				\textup{\text{(\ref{conjecture: maximal collections of alpha[x] are M[x]})}}
			\end{array}
			\arrow[r,Rightarrow, ""']&
			\begin{array}{c}
				\textup{\textbf{\hyperlink{(C2)}{(C2)}}}
			\end{array}
			\arrow[rrr,Leftrightarrow,"\text{Cor. }(\ref{question two is equivalent to the Anderson-Badawi conjectures})"']  &[0em] &[-1.5em] &[0em]
			\begin{array}{c}
				\textup{\text{(\ref{conjecture: graded ring formulation of Anderson-Badawi})}}
			\end{array}
		\end{tikzcd}
	\end{center}
	
	\subsection{Collections of $\alpha[X]$}\label{subsection: collections of alpha[x]}
	
	This subsection focuses on exploring the strictness of the inclusion $\tilde C_j\seq c(\tilde C_j)[X]$ for a collection $C$ of $\alpha[X]$, eventually culminating in the introduction of Conjecture (\ref{conjecture: maximal collections of alpha[x] are M[x]}). 
	This conjecture hypothesizes that this inclusion is an equality whenever $\alpha$ is an ideal of a Noetherian ring and $C$ is a maximal collection of $\alpha$.
	
	We will start this subsection by providing a slight generalization of a well-known theorem by McCoy \cite{mccoy1942remarks}, which states that if $\gamma$ is an ideal of a polynomial ring $A[X]$ such that there exists a nonzero $f\in A[X]$ with $f\gamma=\<0\>$, then there exists a nonzero $a\in A$ such that $a\gamma=\<0\>$. 
	Before we generalize this, we will need to formalize some terminology.
	For a nonzero polynomial $f=\sum_v f_v \prod_{i=1}^n x_i^{v_i}\in A[x_1,...,x_n]=A[X]$, define the degree of $f$ to be the lexicographically largest $v$ such that $f_v\ne 0$, and we define the content ideal $c(f)=\<\{f_v \}_v\>$ to be the ideal of $A$ generated by the coefficients  of $f$. 
	Similarly, if $\beta$ is an ideal of $A[X]$ then we set $c(\beta)=\sum_{f\in \beta}c(f)$.
	
	\begin{prop}\label{Prop: Generalziation of McCoys zero divisor theorem}
		Suppose $\gamma$ is a nonzero ideal of $A[X]=A[x_1,...,x_n]$ and $f\in A[X]$ is nonzero such that $f\gamma=\<0\>$.
		Then there exists some $a\in A$ such that $af\ne 0$ but $c(af)\gamma=\<0\>$.
		\begin{proof}
			We will generalize the proof given in \cite{scott1954divisors}.
			Let $g=\sum_r g_r X^r\in Af$ be nonzero and of (lexicographically) smallest degree possible. 
			With intent to contradict, assume $c(g)\gamma\ne\<0\>$.
			Since $c(g)\gamma\ne \<0\>$ there must exist some element $h=\sum_{v} h_v X^v$ of $\gamma$ such that $c(g)h\ne \<0\>$; equivalently, $c(h)g\ne \<0\>$.
			Let $I$ be the set of indices $v\in I$ with $ h_v\cdot g\ne0$. Notice $gh\in f\gamma=\<0\>$ and so
			\[0=gh=g(\sum_{v\in I}h_v X^v+\sum_{w\not\in I}h_w X^w)=g\sum_{v\in I}h_v X^v=(\sum_r g_r X^r)(\sum_{v\in I}h_v X^v)\]
			Now let $k$ be the largest element of $I$, and let $m$ be largest such that $g_m\ne 0$ (i.e. the leading coefficient of $g$).
			Well, $(\sum_r g_r X^r)(\sum_{v\in I}h_v X^v)=0$ implies $g_m h_k=0$, and thus the degree of $h_kg$ is less than the degree of $g$.
			Since $k\in I$ we know $h_kg\ne 0$, and thus $h_kg$ is a nonzero element of $Af$ whose degree is lexicographically less than the degree of $g$---a contradiction by choice of $g$.
		\end{proof}
	\end{prop}
	
	Now let us apply the above proposition to the setting of collections, which helps us understand the strictness of $\tilde C_j\seq c(\tilde C_j)[X]$.
	\begin{cor}
		Let $C$ be a collection of $\alpha[X]$. For any $f\in \tilde C_j\setminus \alpha[X]$ there exists some $a\in A$ such that $af\not\in \alpha[X]$ but $c(af)\seq \tilde C_j$.
		\begin{proof}
			Apply Proposition (\ref{Prop: Generalziation of McCoys zero divisor theorem}) to the ring $A[X]/\alpha[X]=(A/\alpha)[X]$.
		\end{proof}
	\end{cor}
	
	The following claim can be interpreted as a result in the formate popularized by the Dedekind--Mertens Lemma, but in the context of collections.
	Moreover, this is the strongest general statement we can currently make about the strictness of the inclusion $\tilde C_j\seq c(\tilde C_j)[X]$.
	
	\begin{claim}\label{Claim: c(c_i)^omega(alpha) is contained in tilde c_i}
		Let $\alpha$ be an ideal of a ring $A$ with $\omega(\alpha)<\infty$. For any collection $C$ of $\alpha[X]$ and every $C_j\in C$, we have: \[c\left(\tilde C_j\right)^{\omega(\alpha)}\seq \tilde C_j\]
		\begin{proof}
			Clearly $c(\tilde C_j)^{\omega(\alpha)}\seq \tilde C_j$ if and only if $c(J)^{\omega(\alpha)}\seq \tilde C_j$ for every finitely generated subideal $J$ of $\tilde C_j$. 
			Similarly, $c(J)^{\omega(\alpha)}\seq \tilde C_j=(\alpha[X]:\prod_{i\ne j}C_i)$ if and only if $c(J)^{\omega(\alpha)}\seq (\alpha[X]:J')$ for every finitely generated subideal $J'$ of $\prod_{i\ne j}C_i$.
			We will always have $J J'\seq \alpha[X]$, and thus $c(JJ')\seq \alpha$. 
			
			A weak corollary of the Dedekind--Mertens Lemma states that for all $f,g\in A[X]$ there exists some $k$ such that $c(f)^kc(g)\seq c(fg)$.
			After fixing a finite generating set for both $J=\<\{f_i\}_i\>$ and $J'=\<\{g_j\}_j\>$, we can choose $k$ large enough such that $c(f_i)^kc(g_j)\seq c(f_ig_j)$ for all $i$ and $j$.
			Since $c(f_ig_j)\seq c(JJ')\seq \alpha$, we know $c(f_i)^kc(g_j)\seq \alpha$.
			If $m$ is the size of the finite set $\{f_i\}_i$, then by the pigeonhole principle $c(J)^{mk}=(\sum_i c(f_i))^{km}\seq \sum_i c(f_i)^k$.
			This implies $c(J)^{mk}c(J')\seq \sum_{i,j} c(f_i)^k c(g_j) \seq \alpha$.
			By Theorem (\ref{Theorem: (C1) - the First Anderson-Badawi Conjecture}) we have $c(J)^{\omega(\alpha)}c(J')\seq\alpha$ and therefore $c(J)^{\omega(\alpha)}\seq(\alpha[X]:J')$, as needed.
		\end{proof}
	\end{claim}
	
	Recall that Lemma (\ref{Lemma: omega alpha is finite for Noetherian rings}) states that any ideal of a Noetherian ring will have finite absorbing number, and by Lemma (\ref{Lemma: Max elements of max collections are prime}) the maximal elements of a maximal collection will always be prime ideals.
	Therefore, by the above Claim (\ref{Claim: c(c_i)^omega(alpha) is contained in tilde c_i}), if $\alpha$ is an ideal of a Noetherian ring $A$ then the maximal elements of any maximal collection $M$ of $\alpha$ will be the form $\p[x]$ for some prime ideal $\p$ of $A$.
	With the above results in mind, it becomes natural to ask if \textit{every} element of $M$ will be of the form $\beta[x]$ for some ideal $\beta$ of $A$. Equivalently:
	\begin{conjecture}\label{conjecture: maximal collections of alpha[x] are M[x]}
		For an ideal $\alpha$ of a Noetherian ring, the maximal collections of $\alpha[x]$ are of the form $M[x]:=[M_i[x]]_i$ for some maximal collection $M=[M_i]_i$ of $\alpha$.
	\end{conjecture}
	
	In the Noetherian context, Conjecture (\ref{conjecture: maximal collections of alpha[x] are M[x]}) naturally implies that every maximal collection of $\alpha[X]=\alpha[x_1,...,x_n]$ will be of the form $M[X]$.
	We may also ask if Conjecture (\ref{conjecture: maximal collections of alpha[x] are M[x]}) also holds in the non-Noetherian context, but that will not be explored here.
	
	If $\alpha$ is an ideal of $A$ such that $\omega(\alpha)$ is finite, then we can always find a finitely generated subring $B$ of $A$ large enough such that $\omega(\alpha\cap B)=\omega(\alpha)$.
	Now if $\omega(\alpha[x])> \omega(\alpha)$, then we can choose $B$ large enough such that we also have $\omega((\alpha\cap B)[x])>\omega(\alpha\cap B)$. Therefore, if \hyperlink{(C2)}{(C2)} holds for all finitely generated rings then \hyperlink{(C2)}{(C2)} must hold for all rings. 
	Since finitely generated rings are Noetherian, Conjecture (\ref{conjecture: maximal collections of alpha[x] are M[x]}) naturally implies \hyperlink{(C2)}{(C2)}. 
	
	\subsection{Finitely Many Maximal Collections?}\label{subsection: finitely many maximal collections}
	A primary obstacle one faces when studying absorbing ideals is the difficulty of calculating $\omega(\alpha)$. 
	Currently, there are no known algorithms to compute $\omega(\alpha)$ for any ideal $\alpha$ of $\Z[x_1,...,x_r]$, and we do not even have a proof that such an algorithm can exist.
	
	If an ideal $\alpha$ of a Noetherian ring has finitely many maximal collections, and if there exists some algorithm that computes all its maximal collections, then we would have an algorithm to compute $\omega(\alpha)$.
	There is evidence that such ideals $\alpha$ may be common. 
	For instance, $2$-absorbing ideals of Noetherian rings have finitely many maximal collections and computing them is straightforward:
	By Lemma (\ref{Lemma: Max elements of max collections are prime}) every (non-prime) $2$-absorbing ideal $\alpha$ of a Noetherian ring will have maximal collections of the form $[(\alpha:\p),\p]$ for some associated prime $\p$ (up to reordering). 
	Since any ideal of a Noetherian ring will have finitely many associated primes, we know $2$-absorbing ideals of Noetherian rings have finitely many maximal collections.
	
	With this motivation in mind, we pose the following conjecture in hopes that a positive answer could help motivate a computational study of absorbing ideals and maximal collections.
	\begin{conjecture}\label{conjecture: finitely many maximal collections}
		Any ideal of a Noetherian ring has finitely many maximal collections.
	\end{conjecture}
	To summarize the above discussion, if Conjecture (\ref{conjecture: finitely many maximal collections}) is true then we can ask the natural question: \enquote{\textit{if $\alpha$ is an ideal of a $\Z[x_1,...,x_r]$, does there exist an algorithm which computes all maximal collections of $\alpha$?}}
	
	It turns out that Conjecture (\ref{conjecture: finitely many maximal collections}) is quite strong. 
	Not only does it trivially imply results such as Lemma (\ref{Lemma: omega alpha is finite for Noetherian rings}), but Conjecture (\ref{conjecture: finitely many maximal collections}) is actually stronger than Conjecture (\ref{conjecture: maximal collections of alpha[x] are M[x]}), and is therefore stronger than Anderson and Badawi's second conjecture \hyperlink{(C2)}{(C2)}.
	Before we prove this, we will need some definitions and a lemma.
	
	For a polynomial ring $A[x,y]$, the $x$-degree of an element $f\in A[x,y]$ is the degree of $f=\sum_i f_i x^i$ with coefficients $f_i\in A[y]$. 
	Now let $\beta$ be an ideal of $A[x,y]$, and define the $x$-content $c_x(\beta)$ of $\beta$ to be the content ideal of $\beta$ where $A[x,y]=(A[y])[x]$ is viewed as a polynomial ring in $x$. Hence $\beta\seq c_x(\beta)[x]\seq c(\beta)[x,y]$, and if we define $c_y$ similarly we get $c_y(c_x(\beta)[x])[y]=c(\beta)[x,y]=c_x(c_y(\beta)[y])[x]$.
	For each $h\in A[x]$ define the evaluation map $\mathrm{eval}_h:A[x,y]=(A[x])[y]\rightarrow A[x]$ to be the homomorphism which fixes $A[x]$ and sends $y$ to $h$ (i.e. $\mathrm{eval}_h(\sum_i f_i y^i)=\sum_i f_i h^i$). 
	
	\begin{lemma}\label{Lemma: adjugate matrix cor}
		Let $\alpha$ be an ideal of $A$, and let $\beta$ be an ideal of $A[x,y]$ which can be generated by elements of $x$-degree less than $k$. Then  $\beta\seq \alpha[x,y]$ if and only if $\mathrm{eval}_{x^{k }}(\beta)\seq \alpha[x]$.
		\begin{proof} 
			Let $\{f_i\}_i$ be a generating set of $\beta$ such that the $x$-degree of $f_i$ is less than $k$ for all $f_i$.
			Now let $f_i = \sum_j f_{ij}y^j$ where $f_{ij}\in A[x]$, and thus each $f_{ij}$ must be of $x$-degree less than $k$.
			By Lemma (\ref{equation: key fact for splitting content of polynomials}), if $h\in A[x]$ is of $x$-degree less than $\ell$, then $c(g+x^\ell h)=c(g)+c(h)$ for any $h\in A[x]$.
			Inductively applying this statement to $c(\sum_j f_{ij}x^{kj})$ yields the equality $ c(\sum_j f_{ij}x^{kj})=\sum_j c(f_{ij})$.
			This implies 
			\(c(\mathrm{eval}_{x^{k}}(f_i)) = \sum_jc(f_{ij}) = c(f_i)\), and thus $c(\mathrm{eval}_{x^{k}}(\beta))=c(\beta)$.
			Since $\gamma\seq \alpha[X]$ if and only if $c(\gamma)\seq \alpha$ for any ideal $\gamma$ and set of indeterminant(s) $X$, the equality $c(\mathrm{eval}_{x^{k}}(\beta))=c(\beta)$ implies the result.
		\end{proof}
	\end{lemma}
	
	Recall that for a collection $C$ of $\alpha$, we define $C[x]:=[C_i[x]]_i$, which is always a collection of $\alpha[x]$. 
	For a collection $N$ of $\alpha[x,y]$, we define the content of $N$ to be the list $c(N):=[c(N_i)]_i$ where the content ideals $c(N_i)$ are ideals of $A$---although $c(N)$ need not be a collection of $\alpha$.
	Similarly, we define both $c_x(N)$ and $c_y(N)$ in the same way.
	
	\begin{prop}\label{prop: conjecture finitely many maximal collections implies conjecture maximal collections of alpha[x] are M[x]}
		If $\alpha$ is an ideal of a Noetherian ring $A$ such that $\alpha[x]$ has finitely many maximal collections, then the maximal collections of $\alpha[x]$ will be of the form $M[x]$ for some maximal collection $M$ of $\alpha$. 
		In particular,
		\begin{center}
			Conjecture (\ref{conjecture: finitely many maximal collections}) $\implies$ Conjecture (\ref{conjecture: maximal collections of alpha[x] are M[x]}).
		\end{center}
		\begin{proof}
			First, let us show that it is enough to prove that every maximal collection $M$ of $\alpha[x,y]$ satisfies $M= c_y(M)[y]:=[c_y(M_i)[x]]_i$.
			By symmetry of $x$ and $y$, this implies $M=c_x(M)[x]$. 
			In particular,
			\(c(M_i)[x,y]=c_x(c_y(M_i)[y])[x]= c_x(M_i)[x]=M_i\)
			and thus every maximal collection of $\alpha[x,y]$ will be of the form $M=c(M)[x,y]$. 
			Next, we need to show that this implies all maximal collections $N$ of $\alpha[x]$ are of the form $C[x]$ for some (necessarily maximal) collection $C$ of $\alpha$. 
			Suppose $M$ is a maximal collection of $\alpha[x,y]$ lying above $N[y]$.
			By assumption, we have $M=c(M)[x,y]$ and therefore $c(M)[x]$ is a collection of $\alpha[x]$ lying above $N$. 
			Since $N$ is maximal, we know $N=c(M)[x]$, as needed.
			
			Let $M$ be a maximal collection of $\alpha[x,y]$, we will show $M=c_y(M)[y]$.
			Let $k$ be a positive integer large enough such that every $M_j$ and $\prod_{i\ne j}M_i$ can be generated by elements of $x$-degree less than $\frac{1}{2}k$, which is possible since $A[x,y]$ is Noetherian.
			Since we have assumed $\alpha[x]$ will have finitely many maximal collections, we can pick $k$ large enough such that for every maximal collection $N$ of $\alpha[x]$, each of the ideals $\prod_{i\ne j}N_i$ can be generated by elements of $x$-degree and less than $\frac{1}{2}k$. 
			
			Clearly $\prod_i\mathrm{eval}_{x^{k}}(M_i)\seq\mathrm{eval}_{x^{k}}(\alpha[x,y])= \alpha[x]$, and by Lemma (\ref{Lemma: adjugate matrix cor}) we must have $\prod_{i\ne j}\mathrm{eval}_{x^{k}}(M_i) $ $\not\seq \alpha[x]$ for all $j$.
			Therefore, the list $[\mathrm{eval}_{x^{k}}(M_i)]_i$ is a collection of $\alpha[x]$, and so there must exist some maximal collection $N$ of $\alpha[x]$ which lies above $[\mathrm{eval}_{x^{k}}(M_i)]_i$. 
			By reordering, we may assume $\mathrm{eval}_{x^{k}}(M_i)\seq N_i$ for all $M_i$.
			
			Since $N$ is maximal we have $N_j=(\alpha:\prod_{i\ne j} N_i)$, and thus $\mathrm{eval}_{x^{k}}(M_j)\seq N_j$ implies $\mathrm{eval}_{x^{k}}(M_j)\prod_{i\ne j}N_i\seq \alpha[x]$;
			equivalently, $\mathrm{eval}_{x^{k}}(M_j\prod_{i\ne j}N_i)\seq \alpha[x]$.
			By choice of $k$, each $M_j$ and $\prod_{i\ne j}N_i$ can be generated by elements which are of $x$-degree less than $\frac{1}{2}k$, and therefore there must exist a generating set $\{f_i\}_i$ for the ideal $M_j \prod_{i\ne j}N_i$ such that each $f_i$ has $x$-degree less than $k$. 
			By Lemma (\ref{Lemma: adjugate matrix cor}), we know $M_j\prod_{i\ne j}N_i\seq \alpha[x,y]$. 
			This implies $M_j\seq (\alpha[x,y]:\prod_{i\ne j}N_i)=(\alpha[x]:\prod_{i\ne j}N_i)[y]=N_j[y]$, and thus $M=N[y]$ by maximality of $M$.
			Therefore $M=c_y(M)[y]$, as needed.
		\end{proof}
	\end{prop}
	
	We should observe that \enquote{maximal} can be replaced with \enquote{semimaximal} and the above proof still holds true. 
	Yet there can exist semimaximal collections of some ideal $\alpha[x]$ that are not of the form $N[x]$.
	Indeed, consider the following example influenced by \cite[Proof of Theorem 1]{gilmer1967some}. 
	The list $[\<z-yx\>,\<z+yx\>]$ is a collection of $\<z^2,y^2\>[x]$ in the ring $(\mathbb{Q}[y,z])[x]$, yet the product of the $x$-content $\<z,y\>\cdot\<z,y\>=\<z^2,zy,y^2\>$ is not contained in $\<z^2,y^2\>$.
	Therefore, any semimaximal $2$-collection of $\<z^2,y^2\>[x]$ which lies above $[\<z-yx\>,\<z+yx\>]$ cannot be of the form $N[x]$.
	In particular, this implies $\<z^2,y^2\>[x]$ has infinitely many semimaximal collections.
	Yet each semimaximal collection of $\<z^2,y^2\>[x]$ lies below the unique maximal collection $[\<z,y\>[x],\<z,y\>[x],\<z,y\>[x]]$.
	
	It is rare for an ideal to have finitely many semimaximal collections. 
	In fact, $\alpha$ has finitely many semimaximal collections if and only if there exists finitely many ideals $\beta$ satisfying $\beta=(\alpha:(\alpha:\beta))$. 
	Indeed, for such an ideal $\beta\ne \alpha$ the list $[\beta,(\alpha:\beta)]$ is a semimaximal 2-collection, so the \enquote{only if} direction is trivial.
	It is easily checked that $N_j=(\alpha:(\alpha:N_j))$ for any element $N_j$ of a semimaximal collection $N$. 
	For any finite set of ideals $\mathcal{B}$, the set of collections of $\alpha$ whose elements belong to $\mathcal{B}$ must also be finite. 
	Hence if there are finitely many ideals satisfying $\beta=(\alpha:(\alpha:\beta))$, then there are finitely many semimaximal collections of $\alpha$.
	
	Before we proceed to the next topic, let us sketch a way in which proof of (\ref{prop: conjecture finitely many maximal collections implies conjecture maximal collections of alpha[x] are M[x]}) can be used to simplify \hyperlink{(C2)}{(C2)}.
	First, let (\textdagger) represent the assumption:
	\begin{center}\enquote{\textit{whenever $\alpha$ is an ideal of a Noetherian ring with $\omega(\alpha[x])>\omega(\alpha)$, there are finitely many maximal $\omega(\alpha[x])$-collections of $\alpha[x]$.}}\end{center}
	We will give a sketch of how (\textdagger) is equivalent to \hyperlink{(C2)}{(C2)}, but before we continue we need a result by Laradji---although we give a different proof here.
	\begin{lemma}\textup{\cite[Corollary 2.13]{Laradji}.}\label{lemma: Laradji's Corollary 2.13}
		$\omega(\alpha[x])=\omega(\alpha[X])$ for any ideal $\alpha$.
		\begin{proof}
			By induction, it is enough to show $\omega(\alpha[x])=\omega(\alpha[x,y])$.
			Suppose $\prod_i f_i\in \alpha[x,y]$ defines a collection of $\alpha[x,y]$, and let $k$ be larger than the $x$-degree of $\prod_i f_i$.
			Clearly $\prod_i \mathrm{eval}_{x^k}(f_i)=\mathrm{eval}_{x^k}(\prod_i f_i)\in \alpha[x]$, and by Lemma (\ref{Lemma: adjugate matrix cor}) we know $\prod_{i\ne j} \mathrm{eval}_{x^k}(f_i)  \not\in \alpha[x]$ for all $j$.
			Thus $[\<\mathrm{eval}_{x^k}(f_i)\>]_i$ is a collection of $\alpha[x]$, so $\omega(\alpha[x,y])\leq\omega(\alpha[x])$ and therefore $\omega(\alpha[x])=\omega(\alpha[x,y])$.
		\end{proof}
	\end{lemma}
	Assume (\textdagger) and follow the proof of Proposition (\ref{prop: conjecture finitely many maximal collections implies conjecture maximal collections of alpha[x] are M[x]}) by instead considering collections $M$ of size $\omega(\alpha[x,y])=\omega(\alpha[x])>\omega(\alpha)$ and the finite set of maximal collections $N$ of $\alpha[x]$ with length $\omega(\alpha[x])$. 
	Continuing through the proof of (\ref{prop: conjecture finitely many maximal collections implies conjecture maximal collections of alpha[x] are M[x]}), mutatis mutandis, we get that $M=N[x,y]$ for some collection $N$ of $\alpha$. 
	This is a contradiction since $|N|=|N[x,y]|=|M|>\omega(\alpha)$, and thus $\omega(\alpha[x,y])=\omega(\alpha)$. 
	So (\textdagger) implies \hyperlink{(C2)}{(C2)}, and therefore (\textdagger) is equivalent to \hyperlink{(C2)}{(C2)}.
	
	\subsection{Absorbing Ideals of Graded Rings}\label{subsection: absorbing ideals in graded rings}
	All graded rings are assumed to be graded with a torsion-free abelian group.
	Let $R$ be a ring with $G$-grading $R=\bigoplus_{g\in G} R_g$. 
	An element $y\in R$ is called \textit{homogeneous} if $y\in R_g$ for some $g\in G$, in which case we say $g$ is the (homogeneous) degree of $y$. 
	For $f\in R$ we say $f=\sum_j f_j$ is the \textit{homogeneous decomposition} of $f$ if each $f_j$ is homogeneous of degree $g_j\in G$ with $g_j\ne g_{j'}$ for all $j\ne j'$. 
	We call the (necessary finite) set $\{f_j\}_j$ the \textit{homogeneous components} of $f$.
	
	An ideal $\beta$ of $R$ is called a \textit{homogeneous ideal} if whenever $f\in\beta$ all its homogeneous components $f_j$ must also be in $\beta$. 
	This is equivalent to the condition that $\beta$ can be generated by a set of homogeneous elements.
	For an arbitrary ideal $\alpha$ of $R$, define $\alpha^*$ to be the ideal of $R$ generated by the homogeneous components of $\alpha$, which is the unique largest homogeneous ideal contained inside $\alpha$. i.e.
	\[\alpha^* = \big\<\,\alpha\cap \bigcup_{g\in G} R_g\,\big\>\]
	
	It is well known that if $\p$ is a prime ideal then $\p^*$ must also be a prime ideal \cite[\href{https://stacks.math.columbia.edu/tag/00JT}{Lemma 00JT}]{stacks-project}. 
	Since absorbing ideals are a generalization of prime ideals, it is natural to ask if the same is true of $n$-absorbing ideals:
	
	\begin{conjecture}\label{conjecture: graded ring formulation of Anderson-Badawi}
		For any ideal $\alpha$ of a graded ring $R$, we have $\omega(\alpha^*)\leq \omega(\alpha)$.
	\end{conjecture}
	
	Soon we will see that Conjecture (\ref{conjecture: graded ring formulation of Anderson-Badawi}) is equivalent to Anderson and Badawi's second conjecture \hyperlink{(C2)}{(C2)}. 
	First, let us establish the strongest bound for $\omega(\alpha^*)$ that we can currently prove.
	
	\begin{prop}\label{Lemma: w(alpha*) leq w(alpha[x])}
		Let $\alpha$ be an ideal of a graded ring $R$. Then, \[\omega(\alpha^*)\leq \omega(\alpha[x])\]
		\begin{proof}		
			Let $n=\omega(\alpha[x])$, assume $n<\infty$, and
			suppose $\prod_{i=1}^{n+1} f_i \seq \alpha^*$. We will show that $\prod_{i\ne t} f_i \seq \alpha^*$ for some $t$, and thus $\omega(\alpha^*)\leq n$.
			
			Let $f_i = \sum_{j} f_{ij}$ be the homogeneous decomposition of $f_i$ for each $i$. 
			Let $B$ be the finitely generated subring of $R$ generated by the elements $f_{ij}$ for each $i$ and $j$, and let $H$ be the subgroup generated by the degrees of each $f_{ij}$.
			Clearly $B$ has an $H$-grading induced by the grading on $R$, and each $f_i$ must be an element of $B$.
			This implies that $(\alpha\cap B)^* = \alpha^*\cap B$.

			Since $G$ is torsion-free, we know $H$ must be a finitely generated torsion-free abelian group. By the fundamental theorem of finitely generated abelian groups, we know $H\cong (\Z^r,+)$ for some $r$. Let $g_1,...,g_r$ be a set of (free) generators for $H$, and let $\pi:H\rightarrow \Z^r$ be the isomorphism given by $\pi(g_i^k) = (0,..,k,..,0)$ where $k$ is in the $i^{th}$ spot.
			Using $\pi$, we can create an inclusion $\phi:B\rightarrow B[X^{\pm 1}]:=B[x_1^{\pm 1},...,x_r^{\pm 1}]$ given by 
			\[\phi(y)=\sum_j y_j \prod_k x_k^{\pi(h_j)_k}\]
			where $y=\sum_j y_j$ is the homogeneous decompositions of $y$ with $y_j$ being of degree $h_j\in H$.
			Notice that $\phi$ is an injective homomorphism of graded rings satisfying $\<\phi(\beta)\cap \beta[X^{\pm1}]\>=\<\phi(\beta^*)\>=\beta^*[X^{\pm1}]$ for any ideal $\beta$ of $B$. 
			In particular, \(\phi(\alpha^*\cap B)\seq (\alpha^*\cap B)[X^{\pm1}]\seq (\alpha\cap B)[X^{\pm1}]\).
			
			By Lemma (\ref{lemma: Laradji's Corollary 2.13}) we know that $\omega\left((\alpha\cap B)[x]\right) = \omega\left((\alpha\cap B)[X]\right)$, and so $\omega\left((\alpha\cap B)[X]\right) \leq \omega(\alpha[x])=n$.
			Localization can only decrease the absorbing number, and thus we have $\omega\left((\alpha\cap B)[X^{\pm1}]\right) \leq n$.
			Yet $\phi(\prod_{i=1}^{n+1}f_i)=\prod_{i=1}^{n+1}\phi(f_i) \seq \phi(\alpha^*\cap B)\seq (\alpha\cap B)[X^{\pm1}]$ and so there must exist some $t$ such that 
			\[\phi(\prod_{i\ne t}f_i)=\prod_{i\ne t}\phi(f_i)\seq(\alpha\cap B)[X^{\pm1}]\]
			Therefore, we have $c(\phi(\prod_{i\ne t}f_i))\seq \alpha\cap B$. 
			By construction of $\phi$, given any $y\in B$ we have $c(\phi(y))=\<\{y_j\}_j\>\lhd B$ where $y=\sum_j y_j$ is the homogeneous decomposition of $y$.
			In particular, $c(\phi(\prod_{i\ne t}f_i))\seq \alpha\cap B$ implies every element in the homogeneous decomposition of $\prod_{i\ne t}f_i$ must be an element of $\alpha\cap B$, and subsequently $\prod_{i\ne t}f_i\in (\alpha\cap B)^*\seq\alpha^*$.
			Therefore, $\omega(\alpha^*)\leq \omega(\alpha[x])$.
		\end{proof}
	\end{prop}
	
	As a corollary of the above proposition, we get:
	
	\begin{cor} \label{question two is equivalent to the Anderson-Badawi conjectures}
		Conjecture (\ref{conjecture: graded ring formulation of Anderson-Badawi}) is equivalent to \hyperlink{(C2)}{(C2)}.
		\begin{proof}
			By Proposition (\ref{Lemma: w(alpha*) leq w(alpha[x])}) we know that \hyperlink{(C2)}{(C2)} implies Conjecture (\ref{conjecture: graded ring formulation of Anderson-Badawi}).
			Let us show that Conjecture (\ref{conjecture: graded ring formulation of Anderson-Badawi}) implies \hyperlink{(C2)}{(C2)}. Let $\alpha$ be an ideal of a ring $A$, and let $\beta$ be the ideal of $A[x]$ given by $\beta= \alpha[x]+\<x-1\>$. Notice that $A[x]/\beta\cong A/\alpha$, and thus $\omega(\alpha)=\omega(\beta)$. If the polynomial ring $A[x]$ is given the canonical $\Z$-grading, then $\beta^* = \alpha[x]$. Hence,
			\[\omega(\alpha[x])=\omega(\beta^*)\leq \omega(\beta)=\omega(\alpha)\leq \omega(\alpha[x])\] 
			and therefore $\omega(\alpha)=\omega(\alpha[x])$, as needed. 
		\end{proof}
	\end{cor}

	\section*{Acknowledgements}\label{Section: Acknowledgements}
	I am thankful towards Adam Logan for his remarks on a rough draft of this manuscript, and for his attentive support and advice throughout these years.
	I would also like to thank Ayman Badawi for directing me towards useful material, and for his remarks on my earlier (unrelated) approaches to the subject.
	\printbibliography

@article{Choi_Walker,
  title={The radical of an $ n $-absorbing ideal},
  author={Choi, Hyun Seung and Walker, Andrew},
  journal={Journal of Commutative Algebra},
  volume={12},
  number={2},
  pages={171--177},
  year={2020},
  publisher={Rocky Mountain Mathematics Consortium}
}

@article{Choi_Walker_arxiv,
	title={The radical of an n-absorbing ideal},
	author={Choi, Hyun Seung and Walker, Andrew},
	journal={arXiv preprint arXiv:1610.10077},
	year={2016}
}

@incollection{open,
	title={Open problems in commutative ring theory},
	author={Paul-Jean Cahen and Marco Fontana and Sophie Frisch and Sarah Glaz},
	booktitle={Commutative Algebra},
	pages={353--375},
	year={2014},
	publisher={Springer}
}

@article{n-abs,
	title={On n-absorbing ideals of commutative rings},
	author={Anderson, David F. and Badawi, Ayman},
	journal={Communications in Algebra},
	volume={39},
	number={5},
	pages={1646--1672},
	year={2011},
	publisher={Taylor \& Francis}
}

@article{2-abs, title={On 2-absorbing ideals of commutative rings}, volume={75}, number={3}, journal={Bulletin of the Australian Mathematical Society}, publisher={Cambridge University Press}, author={Badawi, Ayman}, year={2007}, pages={417–429}}

@article{Donadze_infinite_fields,
	title={The Anderson-Badawi conjecture for commutative algebras over infinite fields},
	author={Donadze, Guram},
	journal={Indian Journal of Pure and Applied Mathematics},
	volume={47},
	number={4},
	pages={691--696},
	year={2016},
	publisher={Springer}
}

@article{Donadze_radical_formula,
	title={A proof of the Anderson--Badawi {$\text{rad}(I)^n\subseteq I$} formula for  {$n$}-absorbing ideals},
	author={Donadze, Guram},
	journal={Proceedings-Mathematical Sciences},
	volume={128},
	number={1},
	pages={1--6},
	year={2018},
	publisher={Springer}
}

@inproceedings{Laradji,
	title={On n-absorbing rings and ideals},
	author={Laradji, Abdallah},
	booktitle={Colloquium Mathematicum},
	volume={147},
	pages={265--273},
	year={2017},
	organization={Instytut Matematyczny Polskiej Akademii Nauk}
}

@article{Badawi-overview-of-conjectures,
	title={n-Absorbing ideals of commutative rings and recent progress on three conjectures: a survey},
	author={Badawi, Ayman},
	journal={Rings, Polynomials, and Modules},
	pages={33--52},
	year={2017},
	publisher={Springer}
}

@article{McCoy,
	title={A note on finite unions of ideals and subgroups},
	author={McCoy, Neal H.},
	journal={Proceedings of the American Mathematical Society},
	volume={8},
	number={4},
	pages={633--637},
	year={1957},
	publisher={JSTOR}
}

@article{Gottlieb,
	title={On finite unions of ideals and cosets},
	author={Gottlieb, Christian},
	journal={Communications in Algebra},
	volume={22},
	number={8},
	pages={3087--3097},
	year={1994},
	publisher={Taylor \& Francis}
}

@book{matsumura1989commutative,
	title={Commutative ring theory},
	author={Matsumura, Hideyuki},
	year={1989},
	publisher={Cambridge University Press},
	collection={Cambridge Studies in Advanced Mathematics}
}

@misc{stacks-project,
	author       = {The {Stacks project authors}},
	title        = {The Stacks project},
	howpublished = {\url{https://stacks.math.columbia.edu}},
	year         = {2023}
}

@article{nagata1962local,
	title={Local rings.},
	author={Nagata, Masayoshi},
	year={1962},
	publisher={Interscience Publishers, New York}
}

@article{scott1954divisors,
	title={Divisors of zero in polynomial rings},
	author={Scott, W. R.},
	journal={The American Mathematical Monthly},
	volume={61},
	number={5},
	pages={336--336},
	year={1954},
	publisher={JSTOR}
}

@article{gilmer1967some,
	title={Some applications of the Hilfssatz von Dedekind-Mertens},
	author={Gilmer Jr., Robert W.},
	journal={Mathematica Scandinavica},
	pages={240--244},
	year={1967},
	publisher={JSTOR}
}

@article{mccoy1942remarks,
	title={Remarks on divisors of zero},
	author={McCoy, Neil H.},
	journal={The American Mathematical Monthly},
	volume={49},
	number={5},
	pages={286--295},
	year={1942},
	publisher={Taylor \& Francis}
}
\end{document}